\documentclass[10pt, a4paper]{article}
\usepackage{latexsym}
\usepackage{amsmath}
\usepackage{amssymb}
\usepackage[all]{xy}
\usepackage{amsfonts}
\usepackage{verbatim}
\usepackage{amsthm}
\usepackage{mathrsfs}
\usepackage{epsfig}
\usepackage{xy}
\usepackage{array}
\usepackage{stmaryrd}
\usepackage{graphicx,color}
\usepackage{xcolor}
\usepackage{tikz}
\usetikzlibrary{arrows,calc}
\usepackage{etex}
\usepackage{mathdots}
\usepackage{float}
\usepackage{graphics}
\usepackage{pdflscape}

\usetikzlibrary{arrows,decorations.pathmorphing,decorations.pathreplacing,positioning,shapes.geometric,shapes.misc,decorations.markings,decorations.fractals,calc,patterns}

\usepackage{anysize,hyperref}
\input xypic
\xyoption{all}

\usepackage[perpage,symbol]{footmisc}
\usepackage{setspace}
\renewcommand{\cal}{\mathcal}

\newcommand{\Hom}[3]{Hom_{#3}({#1},{#2})}

\newcommand{\rtri}[6]{\xymatrix@C=1.5em{#1\ar[r]^{#4}&#2\ar[r]^{#5}&#3\ar[r]^{#6}&\Sigma #1}}
\newcommand{\tri}[7]{\xymatrix@C=1.5em{#1\ar[r]^{#5}&#2\ar[r]^{#6}&#3\ar[r]^{#7}&#4}}

\def\X{\mathscr{X}}
\def\Y{\mathscr{Y}}
\def\add{\mbox{add}}
\def\Ext{\mbox{Ext}}
\def\Hom{\mbox{Hom}}

\begin{document}
\baselineskip=15pt
\title{\Large{\bf Torsion pairs in finite $2$-Calabi-Yau triangulated categories with maximal rigid objects\footnote{Supported by the NSF of China (Grants No.\;11671221)}}}
\medskip
\author{Huimin Chang \quad\quad\quad Bin Zhu}

\date{ Department of Mathematical Sciences
\\  Tsinghua University
\\
    100084 Beijing, P. R. China
\\
 chm14@mails.tsinghua.edu.cn (Chang); bzhu@math.tsinghua.edu.cn (Zhu)}

\maketitle
\def\blue{\color{blue}}
\def\red{\color{red}}

\newtheorem{theorem}{Theorem}[section]
\newtheorem{lemma}[theorem]{Lemma}
\newtheorem{corollary}[theorem]{Corollary}
\newtheorem{proposition}[theorem]{Proposition}
\newtheorem{conjecture}{Conjecture}
\theoremstyle{definition}
\newtheorem{definition}[theorem]{Definition}
\newtheorem{question}[theorem]{Question}
\newtheorem{remark}[theorem]{Remark}
\newtheorem{remark*}[]{Remark}
\newtheorem{example}[theorem]{Example}
\newtheorem{example*}[]{Example}

\newtheorem{construction}[theorem]{Construction}
\newtheorem{construction*}[]{Construction}

\newtheorem{assumption}[theorem]{Assumption}
\newtheorem{assumption*}[]{Assumption}
\def\s{\stackrel}
\def\Longrightarrow{{\longrightarrow}}
\def\A{\mathcal{A}}
\def\B{\mathcal{B}}
\def\C{\mathcal{C}}
\def\D{\mathcal{D}}
\def\T{\mathcal{T}}
\def\R{\mathcal{R}}
\def\P{\mathcal{P}}
\def\S{\mathcal{S}}
\def\H{\mathcal{H}}
\def\U{\mathscr{U}}
\def\V{\mathscr{V}}
\def\M{\mathscr{M}}
\def\N{\mathcal{N}}
\def\W{\mathscr{W}}
\def\X{\mathscr{X}}
\def\Y{\mathscr{Y}}
\def\Z{\mathcal {Z}}
\def\I{\mathcal {I}}
\def\add{\mbox{add}}
\def\Aut{\mbox{Aut}}
\def\coker{\mbox{coker}}
\def\deg{\mbox{deg}}
\def\diag{\mbox{diag}}
\def\dim{\mbox{dim}}
\def\End{\mbox{End}}
\def\Ext{\mbox{Ext}}
\def\Hom{\mbox{Hom}}
\def\Gr{\mbox{Gr}}
\def\id{\mbox{id}}
\def\Im{\mbox{Im}}
\def\ind{\mbox{ind}}
\def\mod{\mbox{mod}}
\def\mul{\multiput}
\def\c{\circ}
\def \text{\mbox}

\hyphenation{ap-pro-xi-ma-tion}

\baselineskip=17pt
\parindent=0.5cm


\begin{abstract}
\baselineskip=16pt
We give a classification of (co)torsion pairs in finite $2$-Calabi-Yau triangulated categories with maximal rigid objects which are not cluster tilting. These finite $2$-Calabi-Yau triangulated categories are divided into two main classes: one denoted by $\A_{n,t}$ called of type $A$, and the other denoted by $D_{n,t}$ called of type $D$ \cite{bpr}. By using the geometric model of torsion pairs in cluster categories of type $A, $ or type $D$ in \cite{hjr1, hjr3}, we give a geometric description of torsion pairs in $\A_{n,t}$ or $D_{n,t}$ respectively, via defining the periodic Ptolemy diagrams. This allows to count the number of (co)torsion pairs in these categories. Finally, we determine the hearts of (co)torsion pairs in all finite $2$-Calabi-Yau triangulated categories with maximal rigid objects which are not cluster tilting via quivers and relations.  \\[0.4cm]
\medskip
\end{abstract}

\section{Introduction}
The notion of torsion pairs in abelian categories was first introduced by Dickson \cite{d}, and  triangulated version goes back to Iyama and Yoshino \cite{iy}.
It plays an important role in the study of the algebraic structure and geometric structure of triangulated categories, and  unifies the notion of t-structures, co-t-structures, cluster tilting subcategories and maximal rigid subcategories. Torsion pairs are used to construct certain abelian structures inside triangulated categories (as the hearts of torsion pairs) after  Nakaoka's work \cite{n1}. For a given triangulated category, one may ask a question: How many abelian subquotient categories can be constructed from the triangulated category? In general, this is a difficult question, one may think the hearts of $t$-structures in a given triangulated category in the sense of
Beilinson-Bernstein-Deligne\cite{bbd}. To attack this question, one plan is to classify torsion pairs in the given triangulated category, and then to determine the hearts of torsion pairs. The aim of the paper is to give a classification of torsion pairs and to determine their hearts for certain finite $2$-Calabi-Yau triangulated categories.

Classification of torsion pairs (equivalently cotorsion pairs) has been studied by many people recently.  Ng gave a classification of torsion pairs in the cluster category of type $A_{\infty}$ by Ptolemy diagrams of an $\infty$-gon\cite{ng}. Holm, J${\o}$rgensen and Rubey gave a classification of torsion pairs in the cluster category of type $A_{n}$ via Ptolemy diagrams of a regular $(n+3)$-gon \cite[Theorem A]{hjr1}, they also did the same work for the cluster category of type $D_{n}$ by Ptolemy diagrams of a regular $2n$-gon \cite[Theorem 1.1]{hjr3} and for cluster tubes \cite[Theorem 1.1]{hjr2}. Zhang, Zhou and Zhu gave a classification of torsion pairs in the cluster category of a marked surface \cite[Theorem 4.5]{zzz}. Zhou and Zhu gave a construction and a  classification of torsion pairs in any $2$-Calabi-Yau triangulated category with cluster tilting objects\cite[Theorem 4.4]{zz2}.

Cluster categories associated with finite dimensional hereditary algebras \cite{bmrrt} (see also \cite{ccs} for type $A$) and the stable categories of the preprojective algebras $\Lambda$ of Dynkin quivers \cite{gls} have been used for the categorification of cluster algebras. These categories are $2$-Calabi-Yau triangulated categories with an important class of objects called cluster-tilting objects, which are the analogues of clusters in cluster algebras. The cluster-tilting objects are closely related to a class of objects called maximal rigid objets. Indeed, cluster-tilting objects are maximal rigid objects, but the converse is not true in general \cite{bikr,bmv, kz}.  For a $2$-Calabi-Yau triangulated category,  either all maximal rigid objects are cluster tilting, or none of them are \cite[Theorem 2.6]{zz1}.

Triangulated categories with finitely many indecomposable objects (which we call finite triangulated categories) are a special class of locally finite triangulated categories. By Amiot \cite{a} and Burban-Iyama-Keller-Reiten \cite{bikr} (see also \cite{bpr}),  finite $2$-Calabi-Yau triangulated categories with non-zero maximal rigid objects have a classification which depends on wether the maximal rigid objects are cluster tilting or not. Standard finite $2$-Calabi-Yau triangulated categories with non-zero maximal rigid objects which are not cluster tilting are exactly the following orbit categories: \begin{itemize}
\item  \emph{(Type A)}\,\,  $\cal A_{n,t}=D^{b}(\mathbb{K}A_{(2t+1)(n+1)-3})/\tau^{t(n+1)-1}[1]$, where $n\geq 1$ and $t>1$;
\item \emph{(Type D)}\; $\cal D_{n,t}=D^{b}(\mathbb{K}D_{2t(n+1)})/\tau^{(n+1)}\varphi^{n}$, where $n, t\geq 1$, and where $\varphi$ is induced by an automorphism of $D_{2t(n+1)}$ of order \emph{2};
\item  \emph{(Type E)}\; $D^{b}(\mathbb{K}E_{7})/\tau^{2}$ and $D^{b}(\mathbb{K}E_{7})/\tau^{5}$.

\end{itemize}

 Recently Buan-Palu-Reiten classified the algebras arising from these triangulated categories as the endomorphism algebras of maximal rigid objects via mutations of quivers with relations \cite[Table 1, Table 2]{bpr}.

In this paper, we use the geometric models of torsion pairs in cluster categories of type $A_n$ or type $D_n$ in \cite{hjr1} or \cite{hjr3} respectively to define a notion of periodic Ptolemy diagrams. This allows us to give a complete classification of torsion pairs in the categories $\cal A_{n,t},$ $ \cal D_{n,t}$.  From this classification, we count the number of torsion pairs in these categories. We also determine the hearts of these torsion pairs. These results, combining
with results in \cite{zz2} give a complete picture on torsion pairs and their hearts in finite $2$-Calabi-Yau triangulated categories with non-zero maximal rigid objects.

The paper is organized as follows: In Section $2$, some basic definitions and related results are recalled and some conclusions are achieved on torsion pairs. In Section $3$, we give a geometric description of torsion pairs in  $\A_{n,t}$, where $n\geq 1$ and $t>1$, in the first subsection. In the second subsection, we count the number of torsion pairs in these categories. In the final subsection, we use the same approach to count the number of torsion pairs in $\A_{n,1}$, which is a $2$-Calabi-Yau finite triangulated category of type $A$ with cluster tilting objects. In Section $4$, we give a geometric description of torsion pairs in $\cal D_{n,t}$, and count the number of torsion pairs in these categories. In the last section, we determine the hearts of torsion pairs in finite $2$-Calabi-Yau triangulated categories with maximal rigid objects which are not cluster tilting.

\textbf{Notation.} Unless stated otherwise, $\mathbb{K}$ will be an algebraically closed field of
characteristic zero. Our categories will be assumed $\mathbb{K}$-linear, Hom-finite, Krull-Remark-Schmidt additive categories. $\add\,T$ denotes the additive closure of $T$. Any subcategory is assumed to be one closed under finite direct sums and direct summands. Let $\X$ and $\Y$ be subcategories of a triangulated category $\cal C$. $\X\ast\Y$ denotes the  extension subcategory of $\X$ by $\Y$, whose objects are by definition the objects $M$ with the triangle $X\to M\to Y\to X[1]$, where $X\in\X$ and $Y\in\Y$. We say Hom$_{\cal C}(\X,\Y)=0$ if Hom$_{\cal C}(X,Y)=0,$ for $X\in \X, Y\in \Y$. A subcategory $\X$ is called an extension closed subcategory provided that $\X\ast\X \subseteq\X.$
For a subcategory $\cal D$ of $\cal C$,  we denote by $\cal D^\perp$ (resp. $^\perp\cal D$) the subcategory whose objects are $M\in\cal C$ satisfying $\Hom_{\cal C}(\D,M)=0$ (resp. $\Hom_{\cal C}(M,\D)=0$). For the sake of convenience, we write $[1]$ for the shift functor in any triangulated category unless other stated, and Ext$_{\cal C}^1(X, Y)=\Hom (X,Y[1])$.

\section{Preliminaries}
Firstly, we recall some basic notions based on \cite{bmrrt, iy, n2}.

\begin{definition}\label{b0}
Let $\X$ and $\Y$ be subcategories of a triangulated category $\cal C$.
\begin{itemize}
\item[$(1)$] The pair $(\X,\Y)$ is a $torsion$ $pair$ if
$$\Hom_{\cal C}(\X,\Y)=0\text{ and }\cal C=\X\ast\Y\text{.}$$
The subcategory $\cal I=\X\cap\Y[-1]$ is called the $core$ of the torsion pair.

\item[$(2)$] The pair $(\X,\Y)$ is a $cotorsion$ $pair$ if
$$\Ext^1_{\cal C}(\X,\Y)=0\text{ and }\cal C=\X\ast\Y[1]\text{.}$$
Moreover, we call the subcategory $\cal I=\X\cap \Y $ the $core$ of the cotorsion pair.

\item[$(3)$] A $t$-$structure$ $(\X, \Y )$ in $\cal C $ is a torsion pair such that $\X$ is closed under $[1]$ (equivalently
$\Y $ is closed under $[-1]$).

\item[$(4)$] A subcategory $\cal T$ is called $rigid$ if $\Ext^1_{\cal C}(\cal T,\cal T)=0$. $\cal T$ is called $maximal$ $rigid$ if $\cal T$ is maximal with respect to this property, i.e., if $\Ext^1_{\cal C}(\cal T\oplus addM, \cal T\oplus add M)=0$, then $M\in\add\, T$.

$T$ is called a $rigid$ $object$ if add$T$ is rigid. $T$ is $maximal$ $rigid$ if add$T$ is maximal rigid.

\item[$(5)$] A functorially finite subcategory $\cal T$ is called cluster tilting if $\cal T=\{X\in\cal C|\Ext^1_{\cal C}(X,\cal T)=0\}=\{X\in\cal C|\Ext^1_{\cal C}(\cal T, X)=0\}$. An object $T$ is a $cluster$ $tilting$ object if add$T$ is a cluster tilting subcategory.
\end{itemize}
\end{definition}

\begin{remark}
By Definition \ref{b0}, we know that a pair $(\X,\Y)$ is a cotorsion pair if and only if $(\X,\Y[1])$ is a torsion pair.
\end{remark}

\begin{remark}
For a cotorsion pair $(\X,\Y)$ with core $\cal I=\X\cap \Y $, it is easy to see that $(\X,\Y)$ is a $t$-structure if and only if $\cal I=\{0\}$ \cite{zz2};  $\X$  is  a cluster tilting
 subcategory if and only if $\cal I=\X= \Y$.
\end{remark}
The following result can be found in \cite[Proposition 2.3]{iy}.

\begin{lemma}
\begin{itemize}
               \item [(1)]Let $\X$ be a contravariantly finite and extension closed subcategory of a triangulated category $\cal C$. Then $(\X,\X^\perp)$ is a torsion pair.
               \item [(2)] Let $\X$ be a covariantly finite and extension closed subcategory of a triangulated category $\cal C$. Then $({^\perp\X},\X)$ is a torsion pair.
\end{itemize}

\end{lemma}

\begin{definition}
 \begin{itemize}
   \item [1.]  A triangulated category $\cal C$ is called $2$-$Calabi$-$Yau$ (shortly   $2$-CY) provided there is a functorially isomorphism
Hom$_{\cal C}(X, Y)\simeq D$Hom$_{\cal C}(Y, X[2])$, for all $X, Y\in \cal C$, where $D=$Hom$_{\mathbb{K}}(-,\mathbb{K})$.
   \item [2.] \cite{xz} A triangulated category $\cal C$ is $locally$ $finite$ if for any indecomposable object $X$, there exists only a finite number of isomorphism classes of indecomposable objects $Y$ such that $\Hom_{\cal C}(X,Y)\neq0$.
 $\cal C$ is called a $finite$ $triangulated$ $category$ if it contains only finitely many indecomposable objects up to isomorphisms.
 \end{itemize}

\end{definition}

\begin{remark} Any finite $2$-CY triangulated category $\cal C$ contains a maximal rigid object (may be zero). If the maximal rigid objects of a connected finite $2$-CY triangulated category $\cal C$ are zero, then any torsion pair  $(\X, \Y)$ is a $t$-structure. It follows that $\X$, $\Y$ are triangulated subcategories and $\cal C=\X\oplus \Y$. This implies that $(\X,\Y)=(\cal C, 0)$ or $(0,\cal C)$ (see Proposition \ref{p4}). So the finite $2$-CY triangulated category $\cal C$ which we consider in this paper is assumed to contain a nonzero maximal rigid object. If $\cal C$ contains a cluster tilting object, then any maximal rigid object  is cluster tilting $\cite{zz1}$. Finite $2$-CY triangulated categories with non-zero maximal rigid objects are divided into two classes: one with cluster tilting objects, one without cluster tilting objects. For the first class, we can apply results in \cite{zz2} to obtain a classification of torsion pairs. So we are interested in triangulated categories in the second class in this paper. For
these triangulated categories, Amiot gave a classification (see \cite{a,bpr})

\end{remark}
\begin{lemma}\emph{\cite[Proposition 2.2]{bpr}}
The standard, finite 2-Calabi-Yau, triangulated categories with non-zero maximal rigid objects which are not cluster tilting are exactly the orbit categories:
\begin{itemize}
\item  \emph{(Type A)}\,\,  $\A_{n,t}=D^{b}(\mathbb{K}\vec{A}_{(2t+1)(n+1)-3})/\tau^{t(n+1)-1}[1]$, where $n\geq 1$ and $t>1$;
\item \emph{(Type D)}\; $\D_{n,t}=D^{b}(\mathbb{K}\vec{D}_{2t(n+1)})/\tau^{(n+1)}\varphi^{n}$, where $n, t\geq 1$, and where $\varphi$ is induced by an automorphism of $D_{2t(n+1)}$ of order \emph{2};
\item  \emph{(Type E)}\; $D^{b}(\mathbb{K}\vec{E}_{7})/\tau^{2}$ and $D^{b}(\mathbb{K}\vec{E}_{7})/\tau^{5}$.

\end{itemize}
\end{lemma}

These categories depend on parameters $n,t$. We note that when $t=1$, $\A_{n,1}$ is also a finite $2$-CY triangulated category, it has cluster tilting objects \cite{bpr} which we are also interested in.

In the following, we have some conclusions for torsion pairs in orbit triangulated categories. Firstly, we recall the definition of orbit categories\cite{k,g}.
\begin{definition}
Let $\cal D$ be a triangulated category and $F\colon \cal D\rightarrow \cal D$ be an autoequivalence.  The $orbit$ $category$ $\mathcal O_{F}:= \cal D/F$ has the same objects as $\D$ and its morphisms from $X$ to $Y$ are in bijection with $\bigoplus_{i\in\ \mathbb{Z}} \mathrm{Hom}_{\D}(X,F^{i}Y)$.
\end{definition}
\begin{lemma}\label{a1}
Let $\cal D$ be a locally finite triangulated category and $F\colon \cal D\rightarrow \cal D$ be an autoequivalence such that the orbit category $\mathcal O_{F}= \cal D/F$ is a triangulated category and the projection functor $\pi: \cal D\rightarrow \mathcal O_{F}$ is a triangle functor.
If $(\X, \Y)$ is a torsion pair in $\mathcal O_{F}$, then $( \pi^{-1}(\X), \pi^{-1}(\Y))$ is a torsion pair in $\mathcal D$.
\end{lemma}

\proof We first show $\pi^{-1}(\X)$ is closed under extensions. For any $Z\in \pi^{-1}(\X)\ast\pi^{-1}(\X)$, there exists a triangle
$X_{1}\rightarrow Z\rightarrow X_{2}\rightarrow X_1[1]$ with $X_{1},X_2\in \pi^{-1}(\X)$ in $\mathcal D$.
Since $\pi$ is a triangle functor, we have that $\pi(X_{1})\rightarrow \pi(Z)\rightarrow \pi(X_{2})\rightarrow \pi(X_{1}[1])$ is a triangle in $\mathcal O_{F}$.
 Thus $\pi(Z)\in \X*\X\subseteq \X$, that is, $Z\in \pi^{-1}(\X)$.

 Since $\cal D$ is locally finite, any subcategory of $\cal D$ is functorially finite. It follows that $\pi^{-1}(\X)$ is functorially finite in $\D$. Since $(\X, \Y)$ is a torsion pair in $\mathcal O_{F}$, we have
$\pi^{-1}(\Y)=\pi^{-1}(\X^{\perp})=\pi^{-1}(\X)^{\perp}.$  \qed
\vspace{2mm}

\begin{definition}\label{a2}
Let $\X$ and $\Y$ be subcategories of a triangulated category $\cal C$, and $F\colon \cal C\rightarrow \cal C$ be an autoequivalence. The pair $(\X,\Y)$ is called an \emph{$F$-periodic} torsion pair if $(\X,\Y)$ is a torsion pair and $\X$ is $F$-periodic i.e., $F\X=\X$ (equivalently, $\Y$ is $F$-periodic).
\end{definition}

\begin{lemma}\label{a3}
Let $\cal D$ be a locally finite triangulated category and $F\colon \cal D\rightarrow \cal D$ be an autoequivalence such that $\mathcal O_{F}= \cal D/F$ is a triangulated category and the projection functor $\pi: \cal D\rightarrow \mathcal O_{F}$ is a triangle functor.
If $(\X, \Y)$ is an $F$-periodic torsion pair in $\mathcal D$, then $(\pi(\X), \pi(\Y))$ is a torsion pair in $\mathcal O_{F}$.
\end{lemma}

\proof Since $\X$ is $F$-periodic and $(\X, \Y)$ is a torsion pair, we have
$$\mathrm{Hom}_{\mathcal O_{F}}(\pi(\X), \pi(\Y))=\bigoplus \limits_{i\in\mathbb{Z}} \mathrm{Hom}_{\mathcal D}(F^{i}\X,\Y)=0.$$ For any object $Z\in\mathcal O_{F}$, let $Z'$ be an object in its preimage in $\mathcal D$, that is, $Z'\in\cal D=\X\ast\Y$. Then there exists a triangle $X\rightarrow Z'\rightarrow Y\rightarrow X[1]$ in $\mathcal T$ with $X\in\X $ and $Y\in\Y$. Since $\pi$ is a triangle functor, we have that
 $\pi(X)\rightarrow Z\rightarrow \pi(Y)\rightarrow \pi(X[1])$ is a triangle in $\mathcal O_{F}$, i.e., $Z\in \pi(\X)\ast\pi(\Y)$. This proves that $(\pi(\X), \pi(\Y))$ is a torsion pair in $\mathcal O_{F}$.\qed
 \vspace{2mm}

The following theorem gives a one-to-one correspondence between $F$-periodic torsion pairs in $\cal D$ and torsion pairs in $\mathcal O_{F}$.

\begin{theorem}\label{a4}
Let $\cal D$ be a locally finite triangulated category and $F\colon \cal D\rightarrow \cal D$ be an autoequivalence such that $\mathcal O_{F}:= \cal D/F$ is a triangulated category and the projection functor $\pi: \cal D\rightarrow \mathcal O_{F}$ is a triangle functor. Then there is a bijection between the following sets:
\begin{enumerate}
  \item[\emph{(1)}] The set of $F$-periodic torsion pairs in $\mathcal D$;
  \item[\emph{(2)}] The set of torsion pairs in $\mathcal O_{F}$.
\end{enumerate}
\end{theorem}

\proof This follows directly from Lemma \ref{a1} and Lemma \ref{a3}. \qed

\begin{proposition}\label{p4}
Let $\cal C$ be a connected finite $2$-CY triangulated category. Then the $t$-structures of $\cal C$ are trivial, i.e. $(\cal C,0)$ or $(0,\cal C)$.
\end{proposition}

\proof Suppose $(\X,\Y)$ is a $t$-structure in $\C$. From definition we have $\X[1]\subseteq \X, \Y[-1]\subseteq \Y$. It follows that $\X[1]=\X, \Y[1]=\Y$, i.e. $\X,\Y$ are triangulated subcategories of $\cal C$, and $\C=\X\oplus \Y$. Thus $\X=\cal C$ and $\Y=0$ or $\X=0$ and $\Y=\cal C$. \qed

\vspace{2mm}

\section{Classification of torsion pairs in $\A_{n,t}$}
In this section, we give a classification of torsion pairs in finite $2$-CY triangulated categories of type $A$ with non-zero maximal rigid objects. These categories are denoted by $\cal A_{n,t}$ \cite{bpr}. When $t=1$, the categories $\cal A_{n,1}$ have cluster tilting objects; when $t>1$, the categories $\cal A_{n,t}$ have non-zero maximal rigid objects which are not cluster tilting.

\subsection{A geometric description of torsion pairs in $\A_{n,t}$ }

Let $\mathcal C_{A_{N-3}}$ be the cluster category of type $A_{N-3}$, where $N=(2t+1)(n+1)$. By the universal property of orbit categories \cite{k}, also by the proof of Lemma 2.4 in \cite{bpr}, we know that there exists a covering functor $\pi\colon\mathcal C_{A_{N-3}}\rightarrow\cal A_{n,t}$, which is a triangle functor. Write $F=\tau^{t(n+1)}$, then $F: \mathcal C_{A_{N-3}}\rightarrow\mathcal C_{A_{N-3}}$ is an autoequivalence.  Since $\tau^{N-2}=[-2]$ in $D^{b}(\mathbb{K}\vec{A}_{N-3})$ by \cite{k} and $\tau=[1]$ in $\mathcal C_{A_{N-3}}$, $\tau$ is of order $N$ and $\tau^{n+1}$ is of order $2t+1$ in $\mathcal C_{A_{N-3}}$. Moreover, gcd$(t, 2t+1)=1$ implies that the order of $F=\tau^{t(n+1)}$ is $2t+1$, and the groups generated by $F$ and by $\tau^{n+1}$ are the same, i.e. $<F>=<\tau^{n+1}>$. Therefore $\A_{n,t}$ can be seen as the orbit category $\mathcal C_{A_{N-3}}/\tau^{n+1}$, and $\pi$ is a $(2t+1)$-covering functor (see \cite{bpr}). By Theorem \ref{a4}, we have the following consequence.

\begin{corollary}
There is a bijection between the set of $\tau^{n+1}$-periodic torsion pairs in $\mathcal C_{A_{N-3}}$ and the set of torsion pairs in $\cal A_{n,t}$.
\end{corollary}

In the following, we recall the description of Ptolemy diagrams based on \cite{hjr1}, and give a correspondence between subcategories of $\A_{n,t}$ and collections of diagonals of $N$-gon.

Let $P_{n}$ be an $n$-gon, we label the vertices of $P_{n}$ clockwise by $1,2,\ldots n$ consecutively, where $n\geq 4$ is a positive integer.  A $diagonal$ is a set of two non-neighbouring vertices $\{\alpha, \beta\}$. Two diagonals $\{\alpha_{1}, \alpha_{2}\}$ and $\{\beta_{1}, \beta_{2}\}$ $cross$ if their end points are all distinct and come in the order $\alpha_{1}$, $\beta_{1}$, $\alpha_{2}$, $\beta_{2}$ when moving around the polygon in one direction or the other.

\begin{definition} Let $\mathfrak U$ be a set of diagonals in the $n$-gon $P_{n}$.

 1\cite{hjr1}. $\mathfrak U$ is called a $Ptolemy$ $diagram$ if for any two crossing diagonals $\alpha=\left\{\alpha_{1}, \alpha_{2}\right\}$ and $\beta=\left\{\beta_{1}, \beta_{2}\right\}$ in $\mathfrak U$, those of $\left\{\alpha_{1}, \beta_{1}\right\}$, $\left\{\alpha_{1}, \beta_{2}\right\}$, $\left\{\alpha_{2}, \beta_{1}\right\}$, $\left\{\alpha_{2}, \beta_{2}\right\}$ which are diagonals are in $\mathfrak U$ (see figure \ref{1} for an example).

 2. Fix a positive integer $k|n$, and $n=k\ell$ for some integer $\ell$.  $\mathfrak U$ is called a $k$-$periodic$ collection of diagonals of $P_{n}$ if for each diagonal $(i, j)\in\mathfrak U$, all diagonals $(i+kr, j+kr)$ \emph{(}modulo $n$\emph{)} for $1\leq r\leq \ell$ are in $\mathfrak U$ (see figure \ref{9} for an example).

 3. $\mathfrak U$ is a $k$-$periodic$ $Ptolemy$ $diagram$ if it is a Ptolemy diagram and is $k$-periodic.

 \end{definition}

There is a bijection between indecomposable objects of the cluster category $\C_{A_{N-3}}$ and diagonals of $N$-gon $P_{N}$ \cite{ccs}. In the following, we don't distinct indecomposable objects and diagonals. The Auslander-Reiten translation $\tau$ acts on diagonals is rotation by one vertex in counterclockwise.
\[ \text{dim\ Ext}_{\C_{A_{N-3}}}^{1}(a,b)=\left\{
\begin{array}{cc}
1  &\text{if\ a\ and\ b\ cross},\\
0  &\text{otherwise}.
\end{array}
\right.
\]
Since $\C_{A_{N-3}}$ has only finitely many indecomposable objects, any subcategory of $\C_{A_{N-3}}$ closed under direct sums and direct summands is completely determined by the set of indecomposable objects it contains.  Then the bijection between indecomposable objects of $\C_{A_{N-3}}$ and diagonals of $P_{N}$ extends to a bijection between subcategories of $\C_{A_{N-3}}$ and sets of diagonals of $P_{N}$.

\begin{figure}
\centering
\begin{tikzpicture}[scale=0.4]
\draw[very thick] (0,0) circle (4) ;
\foreach \a in {1,2,5,6}\draw[very thick, color=black!40] (-45*\a-0:4) -- (-45*\a-90:4) ;
\foreach \a in {1,5}\draw[very thick, color=black!40] (-45*\a-0:4) -- (-45*\a-135:4) ;

\foreach \x in {1,2,...,8} {
\draw (-45*\x+135:4) node {$\bullet$} ;
} ;
\draw (-45*1+135:4.7) node {$1$} ;
\draw (-45*2+135:4.7) node {$2$} ;
\draw (-45*3+135:5.7) node {$3$} ;
\draw (-45*4+135:5.3) node {$4$} ;
\draw (-45*5+135:4.7) node {$5$} ;
\draw (-45*6+135:4.8) node {$6$} ;
\draw (-45*7+135:5.9) node {$7$} ;
\draw (-45*8+135:4.7) node {$8$} ;
\end{tikzpicture}
\caption{ $n=8$}
\label{1}
\end{figure}

For any subcategory $\X$ in $\A_{n,t}$, the preimage under the covering functor $\pi$ is a $\tau ^{n+1}$-periodic subcategory $\widetilde{\X}=\pi^{-1}(\X)$ in $\C_{A_{N-3}}$. Moreover, the subcategory $\widetilde{\X}$ corresponds to the set of diagonals of $N$-gon $P_{N}$ by the discussion above, we still denote the corresponding set of diagonals by $\widetilde{\X}$. The corresponding set $\widetilde{\X}$ of diagonals is $(n+1)$-periodic. In the rest of this section, we always use $(i,j)$ to represent an indecomposable object of  $\C_{A_{N-3}}$ or a diagonal of $N$-gon $P_{N}$ without confusion, and $[(i,j)]$ to represent the image under the functor $\pi$. As a consequence, we have the following result.

\begin{figure}
\centering
\begin{tikzpicture}[scale=0.4]
\draw[very thick] (0,0) circle (4) ;
\foreach \a in {1,3,5,7,9}\draw[very thick, color=black!40] (-72*\a-0:4) -- (-72*\a-72:4) ;

\foreach \x in {1,2,...,5} {
\draw (-72*\x+72:4) node {$\bullet$} ;
} ;
\draw (-72*1+72:4.7) node {$1$} ;
\draw (-72*2+72:4.7) node {$3$} ;
\draw (-72*3+72:5.7) node {$5$} ;
\draw (-72*4+72:5.3) node {$7$} ;
\draw (-72*5+72:4.7) node {$9$} ;
\end{tikzpicture}
\caption{ $n=10,k=2$}
\label{9}
\end{figure}
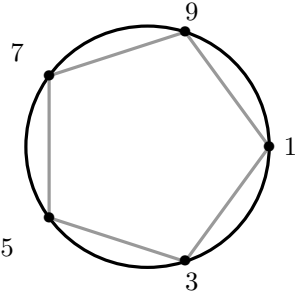

\begin{lemma}
There is a   bijection between the following sets:
\begin{enumerate}
  \item[\emph{(1)}]  Subcategory $\X$ of $\cal A_{n,t}$;
  \item[\emph{(2)}] Collection of diagonals $ \widetilde{\X}$ of the N-gon $P_{N}$ which are $(n+1)$-periodic.
\end{enumerate}
\label{a8}
\end{lemma}

\begin{lemma}\emph{\cite[Theorem A]{hjr1}}
There is a bijection between Ptolemy diagrams of the $(n+3)$-gon and torsion pairs in the cluster category of type $A_{n}$.
\label{a9}
\end{lemma}
The following Lemma gives an equivalent description of torsion pairs in $\cal A_{n,t}$.
\begin{lemma}
Let $\X$ be a subcategory of $\cal A_{n,t}$, and $\widetilde{\X}$ be the corresponding $(n+1)$-periodic collection of diagonals of $N$-gon $P_{N}$. Then the following statements are equivalent:
\begin{enumerate}
  \item[\emph{(1)}]  $(\X, \X^{\bot})$ is a torsion pair in $\cal A_{n,t}$;
  \item[\emph{(2)}] $\X={}^{\bot}(\X^{\bot})$;
  \item[\emph{(3)}] $\widetilde{\X}$ is an $(n+1)$-periodic Ptolemy diagram of N-gon.
\end{enumerate}
\end{lemma}
\proof $(1)\Longleftrightarrow (2)$ is clear.

$"(1)\Rightarrow (3)"$. If $(\X, \X^{\bot})$ is a torsion pair in $\cal A_{n,t}$, then $\X$ corresponds to an $(n+1)$-periodic collection of diagonals $\widetilde{\X}$ of the $N$-gon $P_{N}$ by Lemma \ref{a8}. Moreover, $(\pi^{-1}(\X),\pi^{-1}(\X^{\bot}))$ is a torsion pair in $\C_{A_{N-3}}$ by Lemma \ref{a1}, and $\pi^{-1}(\X)$ corresponds to a Ptolemy diagram of $N$-gon, so $\widetilde{\X}$ is an $(n+1)$-periodic Ptolemy diagram of N-gon.

$"(3)\Rightarrow (1)"$. If $\widetilde{\X}$ is an $(n+1)$-periodic Ptolemy diagram of N-gon, then it corresponds to a  torsion pair in $\C_{A_{N-3}}$ by Lemma \ref{a9}. Moreover, $\widetilde{\X}$ is $(n+1)$-periodic, this implies  the corresponding subcategory $\X$ and $\X^{\bot}$ is a torsion pair in $\A_{n,t}$ by Lemma \ref{a3}.\qed

We will frequently use in this section the coordinate system in the AR-quiver of $\A_{n,t}$, see \cite{bpr} for more details.

\begin{definition}
For a coordinate $(i,j)$ (modulo $N$) with $j>i$ corresponding to an indecomposable object in $\mathcal C_{A_{N-3}}$, we call $j-i-1$ the $level$ of the vertex, and $j-i$ the $length$ of the vertex. Denote the vertex in the AR-quiver of $\cal A_{n,t}$ by $[(i,j)]$ such that all $(i+r(n+1),j+r(n+1))$ for $1\leq r\leq 2t+1$ has to be identified. We also call $j-i-1$ the $level$ of the vertex $[(i,j)]$. The $length$ of the vertex $[(i,j)]$ is $j-i$.
\end{definition}
Buan, Palu and Reiten determined all the indecomposable rigid objects in $\cal A_{n,t}$: for an indecomposable object $[(i,j)]$ in $\A_{n,t}$ with level $j-i-1$, it is rigid if and only if $j-i-1\leq n$ \cite[Lemma 2.4]{bpr}.

\subsection{Torsion pairs in $\A_{n,t}$ with $t>1$}

\begin{proposition}\label{a13}
Let $(\X, \X^{\bot})$ be a torsion pair in $\cal A_{n,t}$, $n\geq 1$, $t>1$, and  $\widetilde{\X}$ be the corresponding $(n+1)$-periodic collection of diagonals of the $N$-gon $P_{N}$. Then precisely one of the following situations occurs:
\begin{enumerate}
  \item[\emph{(1)}] The level of all the indecomposable objects in $\X\leq n$.
  \item[\emph{(2)}] The level of all the indecomposable objects in $\X^{\bot}\leq n$.
\end{enumerate}
\end{proposition}

\proof Note that we can always choose an representative $(i,j)\in[(i,j)]$ such that $1\leq i\leq n+1$, $3\leq j\leq (t+1)(n+1)$ (see Fig.2 in \cite{bpr}). For example, we have $[(1,(t+1)(n+1)+1)]=[(1,1+t(n+1))]$.
\begin{enumerate}

  \item[I.] If the level of all the indecomposable objects of $\X\leq n$,  we claim that $\X^{\bot}$ must contain an element with level$>n$. Indeed, since $\X$ contains only (finitely many) indecomposable rigid objects, we pick an indecomposable object from $\X$ with maximal length. Suppose that its coordinate is $[(1,\ell)]$. Since $(1,\ell)$ corresponds to a rigid object, $3\leq \ell\leq n+2$. If we can show $[(1,\ell+n+1)]\in\X^{\bot}[-1]$, then $\X^{\bot}[-1]$ contains an element with level$>n$, and so does $\X^{\bot}$. Since $t>1$, $\ell+n+1\leq n+2+n+1=2n+3<(t+1)(n+1)$, i.e., $[(1,\ell+n+1)]$ represents a different element from $[(1,\ell)]$ in the AR-quiver of $\cal A_{n,t}$, and the level of $[(1,\ell+n+1)]$ is $\ell+n-1\geq 3+n-1=n+2$. This will complete the proof of our claim. Now we prove that $[(1,\ell+n+1)]\in\X^{\bot}[-1]$. Since $[(1,\ell)]$ is in $\X$ with maximal length, there is no diagonal $(i,j)\in\widetilde{\X}$ with $1<i<\ell<j$, otherwise $(1,\ell)$ has to cross $(i,j)$, but the Ptolemy condition yields a diagonal $(1,j)$ whose length is longer than $(1,\ell)$, a contradiction. If $(1,\ell+n+1)$ crosses a diagonal $(a,b)$ in $\widetilde{\X}$, then $b-a>\ell -1$, a contradiction. This means $(1,\ell+n+1)$ does not cross any diagonal from $\widetilde{\X}$. Similarly, we can prove $(n+2,\ell+2n+2),(2n+3,\ell+3n+3),\ldots \in\widetilde{\X}$, i.e., $[(1,\ell+n+1)]\in\X^{\bot}[-1]$.
  \item[II.] If $\X$ contains an element with level$>n$,  we claim that $\X^{\bot}$ contains only  rigid indecomposable objects. Indeed, without losing generality,
      we suppose $[(1,\ell)]\in\X$ with level $\ell-1-1\geq n+1$, that is, $(t+1)(n+1)\geq \ell\geq n+3$. We choose $\ell$ for different intervals $[n+3,2n+3]$, $[2n+3,3n+4],\ldots,[2t(n+1)+1,(2t+1)(n+1)+1]$.
      \begin{enumerate}
       \item[(a).]  If $\ell\leq 2n+3$, then the corresponding diagonals in $\widetilde{\X}$ are shown
       in figure \ref{2}.

       \item[(b).] If $2n+3\leq\ell\leq3n+4$, then the corresponding diagonals in $\widetilde{\X}$ are shown
       in figure \ref{3}.

       \item[(c).] The other cases are similar.
       \end{enumerate}
       This shows the level of indecomposable objects in $[(1,\ell)]^{\perp}\leq n$, so does $\X^{\perp}$, since $\X^{\perp}\subseteq[(1,\ell)]^{\perp}$.
        \end{enumerate}
      As a consequence, for a torsion pair $(\X, \X^{\bot})$  in $\cal A_{n,t}$, precisely one of $(1)$ and $(2)$ occurs. \qed

\begin{figure}
\centering
\begin{tikzpicture}[scale=0.4]
\draw[very thick] (0,0) circle (4) ;
\foreach \a in {1,3,...,9}
\draw (-36*\a+126:4) edge[very thick, color=black!40, out={-34-36*\a}, in={186-36*\a}] (-36*\a+18:4) ;
\foreach \x in {1,2,...,10} {
\draw (-36*\x+126:4) node {$\bullet$} ;
} ;
\draw (-36*1+126:4.7) node {$1$} ;
\draw (-36*2+126:4.7) node { } ;
\draw (-36*3+126:3) node {\;\;\;\;\;\;\;\;\;\;\;\;\;\;\;\;\;\;$n+2$} ;
\draw (-36*4+126:4.7) node {$\ell$} ;
\draw (-36*5+126:4.9) node {\quad\;$2n+3$} ;
\draw (-36*6+126:4.6) node {$\ell+n+1$} ;
\draw (-36*7+126:4.8) node {$3n+4$\;\;\;} ;
\draw (-36*8+126:6.7) node {$\ell+2n+2$\;\;\;} ;
\draw (-36*9+126:4.7) node { } ;
\draw (-36*10+126:4.7) node { } ;
\end{tikzpicture}
\caption{The maximal length of non-crossing objects is $n+1$}
\label{2}
\end{figure}

\begin{figure}
\centering
\begin{tikzpicture}[scale=0.4]
\draw[very thick] (0,0) circle (4) ;
\foreach \a in {1,3,4,6,9}
\draw (-32.7*\a+122.7:4) edge[very thick, color=black!40, out={-47.3-32.7*\a}, in={172.7-32.7*\a}] (-32.7*\a-8.1:4) ;
\foreach \x in {1,2,...,11} {
\draw (-32.7*\x+122.7:4) node {$\bullet$} ;
} ;
\draw (-32.7*1+122.7:4.6) node {$1$} ;
\draw (-32.7*2+122.7:4.7) node { } ;
\draw (-32.7*3+122.7:4.9) node {\;\;\;$n+2$} ;
\draw (-32.7*4+122.7:5.1) node {\;\;\,$2n+3$} ;
\draw (-32.7*5+122.7:4.5) node {$\ell$} ;
\draw (-32.7*6+122.7:4.6) node {\quad\quad$3n+4$} ;
\draw (-32.7*7+122.7:4.7) node {$\ell+n+1$} ;
\draw (-32.7*8+122.7:5.3) node {$\ell+2n+2$\;\;\;\;\;;\;\;};
\draw (-32.7*9+122.7:4.7) node { } ;
\draw (-32.7*10+122.7:5.4) node {$\ell+3n+3$\;\;\;\;\;\;} ;
\draw (-32.7*11+122.7:4.7) node { } ;
\end{tikzpicture}
\caption{The maximal length of non-crossing objects is $n+1$}
\label{3}
\end{figure}
This Proposition immediately yields the following important conclusion.

\begin{corollary}\emph{\cite{bpr}}
$\cal A_{n,t}$ do not contain any cluster tilting object.

\end{corollary}
\vspace{2mm}

By Proposition \ref{a13}, the classification of torsion pairs $(\X, \X^{\bot})$ in $\cal A_{n,t}$ reduces to the classification of the possible halves $\X$ (or $\X^{\bot}$) of a torsion pair, whose all indecomposable objects are strictly below level $(n+1)$ in the AR-quiver of $\cal A_{n,t}$.

\begin{definition}
Let $(i,j)$ be a diagonal of $N$-gon $P_{N}$. The $wing$ $W(i,j)$ of $(i,j)$ consists of all diagonals $(r,s)$ of the $N$-gon such that $i\leq r\leq s\leq j$, that is all diagonals which are overarched by $(i,j)$.  $[(i,j)]$ represents a vertex in the AR-quiver of $\cal A_{n,t}$, the corresponding wing is denoted by $W[(i,j)]$.
\end{definition}

\begin{theorem}
There are bijections between the following sets:
\begin{enumerate}
  \item[\emph{(1)}] Torsion pairs $(\X, \X^{\bot})$ in $\cal A_{n,t}$ such that the level of all the indecomposable objects in $\X\leq n$;
  \item[\emph{(2)}] $(n+1)$-periodic Ptolemy diagrams $\widetilde{\X}$ of $N$-gon $P_{N}$ such that all diagonals in $\widetilde{\X}$ have length at most $n+1$;
  \item[\emph{(3)}] Collections $\left\{([(i_{1},j_{1})], [W_{1}]), \ldots, ([(i_{r},j_{r})], [W_{r}])\right\}$ of pairs consisting of vertices $[(i_{\ell},j_{\ell})]$ of level$\leq n$ in the AR-quiver of $A_{n,t}$ and subset $[W_{\ell}]\subset W[(i_{\ell},j_{\ell})]$ of their wings such that for any different $k, \ell\in \left\{1,2,\ldots,r\right\}$, we have
  $$ W[(i_{k},j_{k})][1]\cap W[(i_{\ell},j_{\ell})]=\emptyset,$$
  and the $(n+1)$-periodic collection $W_{\ell}$ corresponding to $[W_{\ell}]$ is a Ptolemy diagram.
\end{enumerate}
\label{b}
\end{theorem}
\proof The proof is similar as in the case of cluster tubes \cite[Theorem 4.4]{hjr2}.

Note that the number of indecomposable rigid objects in $\cal A_{n,t}$ is independent of $t$, we have the following result.
\begin{corollary}
The number of torsion pairs in $\cal A_{n,t}$ with $n\geq 1$, $t>1$ is independent of $t$.
\end{corollary}

Therefore counting the number of torsion pairs in $\cal A_{n,t}$ reduces to counting the possible sets of pairs in the AR-quiver of $\cal A_{n,t}:$
 $\left\{([(i_{1},j_{1})], [W_{1}]), \ldots, ([(i_{r},j_{r})], [W_{r}])\right\}$. This is the same as in the process of counting torsion pairs in the cluster tube of rank $n+1$, see \cite{hjr2} for details.
\begin{theorem}
The number of torsion pairs in $\cal A_{n,t}$ with $n\geq 1$, $t>1$ is the same as  the cluster tube of rank $n+1$, that is
$$T_{n+1}=\sum\limits_{\ell\geq 0}2^{\ell+1}\binom{n+\ell}{\ell}\binom{2n+1}{n-2\ell},$$
where $T_{n+1}$ represents the number of torsion pairs in the cluster tube of rank $n+1$.
\label{a}
\end{theorem}

\begin{figure}
\centering
\begin{tikzpicture}[scale=0.5,
fl/.style={->,shorten <=6pt, shorten >=6pt,>=latex}]
\foreach \x in {1,2,3} {
 \pgfmathparse{9-\x}\let\z\pgfmathresult ;
  \foreach \y in {2,3,...,\z} {
   \newcount\u ;
   \pgfmathsetcount{\u}{\x+\y} ;
   \draw (2*\x+\y-4,\y-2) node[scale=0.7] {\x$\;$\the\u} ;
   \draw[fl] (\y-4+2*\x,\y-2) -- (\y-3+2*\x,\y-1) ;
   \draw[fl] (\y-3+2*\x,\y-1) -- (\y-2+2*\x,\y-2) ;
  } ;
} ;
\draw[thick, dashed, blue] (-1.1,-0.5) -- ++(5.6,0) -- ++(4.5,4.5) -- ++(-2.9,2.9) -- cycle ;

\begin{scope}[xshift=8cm]
 \foreach \x in {1,2,3} {
 \pgfmathparse{9-\x}\let\z\pgfmathresult ;
  \foreach \y in {2,3,...,\z} {
   \newcount\w ;
   \pgfmathsetcount{\w}{\x+\y} ;
   \draw (2*\x+\y-6,\y-2) node[scale=0.7] {\x$\;$\the\w} ;
   \draw[fl] (\y-6+2*\x,\y-2) -- (\y-5+2*\x,\y-1) ;
   \draw[fl] (\y-5+2*\x,\y-1) -- (\y-4+2*\x,\y-2) ;
  } ;
} ;
\end{scope}
\begin{scope}[xshift=240]
\draw[thick, dashed, blue] (-3.5,-0.5) -- ++(5.6,0) -- ++(4.5,4.5) -- ++(-2.9,2.9) -- cycle ;
\end{scope}
\begin{scope}[xshift=4cm, yshift=16cm, rotate=180, xscale=-1]
 \foreach \x in {1,2,3} {
 \pgfmathparse{9-\x}\let\z\pgfmathresult ;
  \foreach \y in {2,3,...,\z} {
   \newcount\r ;
   \pgfmathsetcount{\r}{\x+\y} ;
   \draw (2*\x+\y-5,\y+3) node[scale=0.7] {\x$\;$\the\r} ;
   \draw[fl] (\y-5+2*\x,\y+3) -- (\y-4+2*\x,\y+4) ;
   \draw[fl] (\y-4+2*\x,\y+4) -- (\y-3+2*\x,\y+3) ;
  } ;
} ;
\end{scope}
\begin{scope}[xshift=4.2cm, yshift=16.2cm, rotate=180, xscale=-1]
\draw[thick, dashed, blue] (-2.2,+4.8) -- ++(5.5,-0.1) -- ++(4.5,4.5) -- ++(-2.9,2.9) -- cycle ;
\end{scope}
\begin{scope}[xshift=12cm, yshift=16cm, rotate=180, xscale=-1]
 \foreach \x in {1,2,3} {
 \pgfmathparse{9-\x}\let\z\pgfmathresult ;
  \foreach \y in {2,3,...,\z} {
   \newcount\r ;
   \pgfmathsetcount{\r}{\x+\y} ;
   \draw (2*\x+\y-7,\y+3) node[scale=0.7] {\x$\;$\the\r} ;
   \draw[fl] (\y-7+2*\x,\y+3) -- (\y-6+2*\x,\y+4) ;
   \draw[fl] (\y-6+2*\x,\y+4) -- (\y-5+2*\x,\y+3) ;
  } ;
 \draw (3.5+\x,12.5-\x) node[circle, fill=white, scale=1.2] {} ;
} ;
\end{scope}
\begin{scope}[xshift=12.2cm, yshift=16.2cm, rotate=180, xscale=-1]
\draw[thick, dashed, blue] (-4.3,+4.7) -- ++(5.6,0) -- ++(4.5,4.5) -- ++(-2.9,2.9) -- cycle ;
\end{scope}
\draw (3,5) node[scale=1.5] {$\cdots$} ;
\draw (3,9) node[scale=1.5] {$\cdots$} ;
\draw (17,4) node[scale=1.5] {$\cdots$} ;
\end{tikzpicture}
\caption{ The AR-quiver of $\cal A_{2,2}$}
\label{4}
\end{figure}
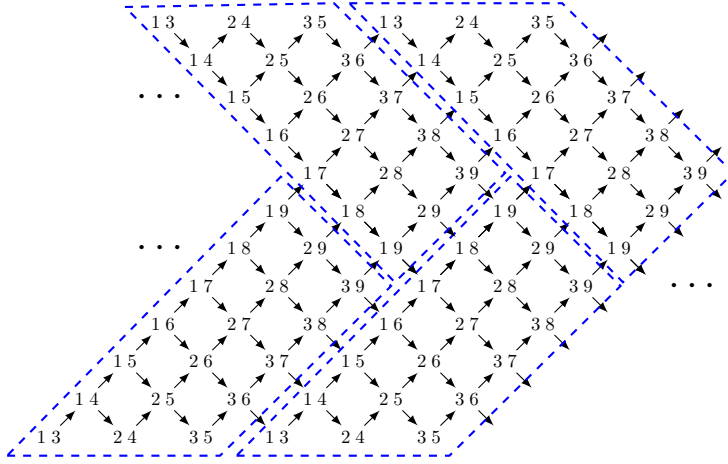

\begin{example}\label{c8}
When
$n=2,t=2$, $\cal A_{2,2}=D^b(\mathbb{K}\vec{A}_{12})/\tau^5[1]$ is  $2$-CY with non-zero maximal rigid objects, whose Auslander-Reiten quiver is shown in figure \ref{4}.  From Theorem \ref{a}, we have $T_{2+1}=32$. Now we construct all torsion pairs by Theorem \ref{b}.
\begin{eqnarray*}
\X_{1}=\{[(0)] \}&\quad\quad& \X^\perp_{1}[-1]=\cal A_{2,2}\nonumber \\
\X_{2}=\{[(13)] \}&\quad\quad& \X^\perp_{2}[-1]=\{[(13)],[(14)],[(16)],[(17)],[(19)],[(36)],[(37)],[(39)]\}\nonumber \\
\X_{3}=\{[(24)] \}&\quad\quad&\X^\perp_{3}[-1]=\{[(14)],[(15)],[(17)],[(18)],[(24)],[(25)],[(27)],[(28)]\}\nonumber \\
\X_{4}=\{[(35)] \}&\quad\quad&\X^\perp_{4}[-1]=\{[(25)],[(26)],[(28)],[(29)],[(35)],[(36)],[(38)],[(39)]\}\nonumber \\
\X_{5}=\{[(14)]\}&\quad\quad& \X^\perp_{5}[-1]=\{[(13)],[(14)],[(17)],[(24)]\}\nonumber \\
\X_{6}=\{[(25)]\}&\quad\quad& \X^\perp_{6}[-1]=\{[(24)],[(25)],[(28)],[(35)]\}\nonumber \\
\X_{7}=\{[(36)]\}&\quad\quad& \X^\perp_{7}[-1]=\{[(13)],[(35)],[(36)],[(39)]\}\nonumber \\
\X_{8}=\{[(13)],[(14)] \}&\quad\quad& \X^\perp_{8}[-1]=\{[(13)],[(14)],[(17)]\}\nonumber \\
\X_{9}=\{[(24)],[(25)] \}&\quad\quad& \X^\perp_{9}[-1]=\{[(24)],[(25)],[(28)]\}\nonumber \\
\X_{10}=\{[(35)],[(36)] \}&\quad\quad& \X^\perp_{10}[-1]=\{[(35)],[(36)],[(39)]\}\nonumber \\
\X_{11}=\{[(14)],[(24)]\}&\quad\quad& \X^\perp_{11}[-1]=\{[(14)],[(17)],[(24)]\}\nonumber \\
\X_{12}=\{[(25)],[(35)] \}&\quad\quad& \X^\perp_{12}[-1]=\{[(25)],[(28)],[(35)]\}\nonumber \\
\X_{13}=\{[(13)],[(36)] \}&\quad\quad& \X^\perp_{13}[-1]=\{[(13)],[(36)],[(39)]\}\nonumber \\
\X_{14}=\{[(13)],[(14)],[(24)] \}&\quad\quad& \X^\perp_{14}[-1]=\{[(14)],[(17)]\}\nonumber \\
\X_{15}=\{[(24)],[(25)],[(35)]\}&\quad\quad& \X^\perp_{15}[-1]=\{[(25)],[(28)]\}\nonumber \\
\X_{16}=\{[(13)],[(35)],[(36)] \}&\quad\quad& \X^\perp_{16}[-1]=\{[(36)],[(39)]\}\nonumber \\
\end{eqnarray*}
\end{example}
\vspace{-0.75cm}
In this example, the collection of pairs in Theorem \ref{b} $(3)$ has one element,i.e., $r=1$. The pair $\{([(i_{1},j_{1})], [W_{1}])\}$ with $[W_{1}]$  containing zero object is the subcategory $\X_{1}$; The pairs $\{([(i_{1},j_{1})], [W_{1}])\}$ with $[W_{1}]$  containing one object are the following subcategories: $\X_{2}$ which corresponds to  $([(1,3)],\{[(1,3)]\})$, $\X_{3}$ which corresponds to  $([(2,4)],\{[(2,4)]\})$, $\X_{4}$ which corresponds to   $([(3,5)],\{[(3,5)]\})$, $\X_{5}$ which corresponds to   $([(1,4)],\{[(1,4)]\})$, $\X_{6}$ which corresponds to  $([(2,5)],\{[(2,5)]\})$, $\X_{7}$ which corresponds to   $([(3,6)],\{[(3,6)]\})$. The pairs $\{([(i_{1},j_{1})], [W_{1}])\}$ with $[W_{1}]$  containing two objects are the subcategories: $\X_{8}$ which corresponds to  $([(1,4)],\{[(1,3)],[(1,4)]\})$, $\X_{9}$ which corresponds to  $([(2,5)],\{[(2,4)],[(2,5)]\})$, $\X_{10}$ which corresponds to $([(3,6)],\{[(3,5)],[(3,6)]\})$, $\X_{11}$ which corresponds to $([(1,4)],\{[(2,4)],[(1,4)]\})$, $\X_{12}$ which corresponds to  $([(2,5)],\{[(2,5)],[(3,5)]\})$, $\X_{13}$ which corresponds to $([(3,6)],\{[(1,3)],[(3,6)]\})$. The pairs $\{([(i_{1},j_{1})], [W_{1}])\}$ with $[W_{1}]$  containing three objects are the subcategories: $\X_{14}$ which corresponds to $([(1,4)],\{[(1,3)],[(1,4)],[(2,4)]\})$, $\X_{15}$ which corresponds to  $([(2,5)],\{[(2,5)],[(3,5)],[(2,4)]\})$, $\X_{16}$ which corresponds to  $([(3,6)],\{[(1,3)],[(3,5)],[(3,6)]\})$.

Then $(\X_i, \X^\perp_{i}[-1])$ is a cotorsion pair in $\cal A_{2,2}$, where $i=1,2,3,\cdots,16$.
It follows that $(\X_i, \X^\perp_{i})$ is a torsion pair in $\cal A_{2,2}$, for any $i$.
Similarly, we can know that $({^\perp\X_{i}},\X_i)$ is a torsion pair in $\cal A_{2,2}$, for $i=1,2,3,\cdots,16$.

\subsection{Torsion pairs in $\A_{n,1}$}
By Amiot \cite{a}, Burban-Iyama-Keller-Reiten \cite{bikr}, the standard, finite $2$-CY triangulated categories of type $A$ with cluster tilting objects are cluster categories of type $A$ and the orbit categories $D^{b}(\vec{A}_{3n})/\tau^{n}[1]$ with $n\geq 1$, which is $\A_{n,1}$ (see Proposition 2.1 in \cite{bpr}). There is also a covering functor from the cluster category $\mathcal C_{A_{3n}}$ to $\cal A_{n,1}$. We have the following result.

\begin{theorem}
The number of torsion pairs in $\cal A_{n,1}$ is
$$N_{n,1}=T_{n+1}-t_{n,1}=\sum\limits_{\ell\geq 0}2^{\ell+1}\binom{n+\ell}{\ell}\binom{2n+1}{n-2\ell}-\sum\limits_{\ell\geq 0}2^{\ell}\binom{n+\ell}{\ell}\binom{2n}{n-2\ell},$$
where $T_{n+1}$ represents the number of torsion pairs in the cluster tube of rank $n+1$, and $t_{n,1}=(n+1)s_{n+2}$, where $s_{n+2}$ represents the number of torsion pairs in the cluster category of type $A_{n-1}$.
\label{d}
\end{theorem}

\proof Recall that an object $[(i,j)]$ in $\cal A_{n,1}$ is rigid if and only if its length $\leq n+1$ \cite{bpr}. For a torsion pair $(\X,  \Y)$ in $\cal A_{n,1}$, if $\X$ (resp. $\Y$) contains a diagonal whose length is longer than $n+1$, i.e. non-rigid object, then $\Y$ (resp. $\X$) contains only indecomposable rigid objects. The proof is the same as II in the proof of Proposition \ref{a13}. Thus, we have two subclasses of torsion pairs in $\cal A_{n,1}$:
\begin{itemize}
  \item [(I)] Torsion pairs $(\X,  \Y)$ such that $\X$ contains only indecomposable rigid objects;
  \item [(II)] Torsion pairs $(\X,  \Y)$ such that $\Y$ contains only indecomposable rigid objects.
\end{itemize}
The intersection of class (I) and class (II) is the subclass of torsion pairs $(\X,  \Y)$ in $\cal A_{n,1}$ such that both $\X$ and $\Y$ contain only indecomposable rigid objects. Next, we consider the case that both $\X$ and $\Y$ contain only indecomposable rigid objects in the following.

If $\X$ is a cluster tilting subcategory in $\cal A_{n,1}$, then $\cal A_{n,1}$ has a cotorsion pair $(\X,  \X)$. Besides these, there are some other torsion pairs $(\X,  \Y)$ with both  $\X$ and $\Y$ containing only indecomposable rigid objects. We give a characterization of them below.

Claim: for a torsion pair $(\X,  \Y)$, both  $\X$ and $\Y$ contain only indecomposable rigid objects if and only if $\X$ and $\Y$ contain one of the rigid objects $[(i,(i+n+1))]$ with $i\in\left\{1,2,\ldots,n+1\right\}$, where these rigid objects are the indecomposable rigid objects with maximal length in $\cal A_{n,1}$. We prove the claim  in the following.
\vspace{2mm}

Let $\widetilde{\X}$ be the corresponding $(n+1)$-periodic collection of diagonals of the $N$-gon with $N=(2t+1)(n+1)=3(n+1)$ as before. If $\X$ contains one  rigid object $[(i,(i+n+1))]$ with maximal length, where $i\in\left\{1,2,\ldots,n+1\right\}$. Without losing generality, we assume $[(1,n+2)]\in\X$ (up to shifting). Then $(1,n+2)$, $(n+2,2n+3)$ and $(2n+3,3n+4)=(2n+3,1)$ are in $\widetilde{\X}$, so the maximal length of diagonals that do not cross with diagonals in $\widetilde{\X}$ is $n+1$. Thus $\Y[1]$ contains only indecomposable rigid objects, so does $\Y$. Similarly one can prove that $\X$ contains only rigid indecomposable objects. Conversely, if $(\X,  \Y)$ is a torsion pair with both  $\X$ and $\Y$ contain only indecomposable rigid objects, we pick a diagonal from $\X$ with maximal length. Suppose its coordinate is $[(1,\ell)]$. It follows from $[(1,\ell)]$ corresponding to an indecomposable rigid object that $3\leq \ell\leq n+2$. If $\ell<n+2$, then $[(1,\ell+n+1)]$ is in $\X^{\bot}[-1]$, which is similar to the part I in the proof of Proposition  \ref{a13} , i.e., $\Y$ contains an indecomposable object with length $>n+1$, which is non-rigid, a contradiction. Thus $l=n+2$, and $\X$ contains the diagonal $[(1,n+2)]$. Same proof implies that $\Y$ contains a rigid object $[(i,(i+n+1))]$, where $i\in\left\{1,2,\ldots,n+1\right\}$.  This completes the proof of the claim above.

It is easy to prove that for a torsion pair $(\X,\Y)$ with $\X$ and $\Y$ containing only rigid indecomposable objects, both $\X$ and $\Y$ contain precisely one indecomposable rigid object with maximal length in $\cal A_{n,1}$. Otherwise, for any such two rigid objects in $\X$, the corresponding diagonals cross, and produce a diagonal with length bigger than $n+1$, which is non-rigid, a contradiction.

Thus the number of torsion pairs in $\A_{n,1}$ is the number of torsion pairs in class (I) plus the number of those in class (II) and minus the number of torsion pairs in the intersection of class (I) and class (II). The sum of numbers of torsion pairs in class (I) and of those in class (II) is the same as the number of torsion pairs in  $\A_{n,t}$ with $t>1$ (compare to Theorem 3.12). We need to count the number of torsion pairs in the intersection. The intersection of class (I) and class (II) is the set of torsion pairs $(\X,\Y)$ with $\X$ and $\Y$ containing only indecomposable rigid objects, which contains one and only one of $[(i,(i+n+1))]$, where $i\in\left\{1,2,\ldots,n+1\right\}$. For the part $\X$ of such a torsion pair $(\X,\Y)$,  the corresponding $(n+1)$-periodic set $\widetilde{\X}$ of diagonals of the $N$-gon contains one diagonal of maximal length. Without losing generality, we assume that it contains $[(1,n+2)]$. Then the diagonals in $\widetilde{\X}$  can be written as $[(i,j)]$ with $1\leq i<j\leq n+2$, which does not cross with $[(1,n+2)]$. Then we have a Ptolemy diagram of $n+2$-gon consisting of diagonals $(i,j)$ such that $[(i,j)] \in \widetilde{\X}\setminus \{[(1,n+2)]\}$, and $1\leq i<j\leq n+2$. Any Ptolemy diagram of $n+2$-gon  gives an $(n+1)$-periodic Ptolemy diagram $\widetilde{\X}$ of the $N$-gon by adding a diagonal $[(1,n+2)].$  Thus, the number of the intersection of class (I) and class (II) is $t_{n,1}=(n+1)s_{n+2}$, where $s_{n+2}$ is the number of Ptolemy diagrams of the $(n+2)$-gon, $s_{n+2}=\frac{1}{n+1}\sum\limits_{\ell\geq 0}2^{\ell}\binom{n+\ell}{\ell}\binom{2n}{n-2\ell}$ by \cite{hjr1}. Then we get the conclusion.  \qed

\begin{figure}
\centering
\begin{tikzpicture}[scale=0.5,
fl/.style={->,shorten <=6pt, shorten >=6pt,>=latex}]
\foreach \x in {1,2,3} {
 \pgfmathparse{6-\x}\let\z\pgfmathresult ;
  \foreach \y in {2,3,...,\z} {
   \newcount\u ;
   \pgfmathsetcount{\u}{\x+\y} ;
   \draw (2*\x+\y-4,\y-2) node[scale=0.7] {\x$\;$\the\u} ;
   \draw[fl] (\y-4+2*\x,\y-2) -- (\y-3+2*\x,\y-1) ;
   \draw[fl] (\y-3+2*\x,\y-1) -- (\y-2+2*\x,\y-2) ;
  } ;
} ;
\draw[thick, dashed, blue] (-1.1,-0.5) -- ++(5.6,0) -- ++(1.5,1.5) -- ++(-2.9,2.9) -- cycle ;

\begin{scope}[xshift=8cm]
 \foreach \x in {1,2,3} {
 \pgfmathparse{6-\x}\let\z\pgfmathresult ;
  \foreach \y in {2,3,...,\z} {
   \newcount\w ;
   \pgfmathsetcount{\w}{\x+\y} ;
   \draw (2*\x+\y-6,\y-2) node[scale=0.7] {\x$\;$\the\w} ;
   \draw[fl] (\y-6+2*\x,\y-2) -- (\y-5+2*\x,\y-1) ;
   \draw[fl] (\y-5+2*\x,\y-1) -- (\y-4+2*\x,\y-2) ;
  } ;
} ;
\end{scope}
\begin{scope}[xshift=240]
\draw[thick, dashed, blue] (-3.5,-0.5) -- ++(5.6,0) -- ++(1.5,1.5) -- ++(-2.9,2.9) -- cycle ;
\end{scope}
\begin{scope}[xshift=4cm, yshift=16cm, rotate=180, xscale=-1]
 \foreach \x in {1,2,3} {
 \pgfmathparse{6-\x}\let\z\pgfmathresult ;
  \foreach \y in {2,3,...,\z} {
   \newcount\r ;
   \pgfmathsetcount{\r}{\x+\y} ;
   \draw (2*\x+\y-5,\y+9) node[scale=0.7] {\x$\;$\the\r} ;
   \draw[fl] (\y-5+2*\x,\y+9) -- (\y-4+2*\x,\y+10) ;
   \draw[fl] (\y-4+2*\x,\y+10) -- (\y-3+2*\x,\y+9) ;
  } ;
} ;
\end{scope}
\begin{scope}[xshift=4.2cm, yshift=16.2cm, rotate=180, xscale=-1]
\draw[thick, dashed, blue] (-2.2,+10.8) -- ++(5.5,-0.1) -- ++(1.5,1.5) -- ++(-2.9,2.9) -- cycle ;
\end{scope}
\begin{scope}[xshift=12cm, yshift=16cm, rotate=180, xscale=-1]
 \foreach \x in {1,2,3} {
 \pgfmathparse{6-\x}\let\z\pgfmathresult ;
  \foreach \y in {2,3,...,\z} {
   \newcount\r ;
   \pgfmathsetcount{\r}{\x+\y} ;
   \draw (2*\x+\y-7,\y+9) node[scale=0.7] {\x$\;$\the\r} ;
   \draw[fl] (\y-7+2*\x,\y+9) -- (\y-6+2*\x,\y+10) ;
   \draw[fl] (\y-6+2*\x,\y+10) -- (\y-5+2*\x,\y+9) ;
  } ;
 \draw (0.5+\x,15.5-\x) node[circle, fill=white, scale=1.2] {} ;
} ;
\end{scope}
\begin{scope}[xshift=12.2cm, yshift=16.2cm, rotate=180, xscale=-1]
\draw[thick, dashed, blue] (-4.4,+10.7) -- ++(5.6,0) -- ++(1.5,1.5) -- ++(-2.8,2.8) -- cycle ;
\end{scope}
\draw (0,4) node[scale=1.5] {$\cdots$} ;
\draw (17,4) node[scale=1.5] {$\cdots$} ;
\end{tikzpicture}
\caption{ The AR-quiver of $\cal A_{2,1}$}
\label{5}
\end{figure}
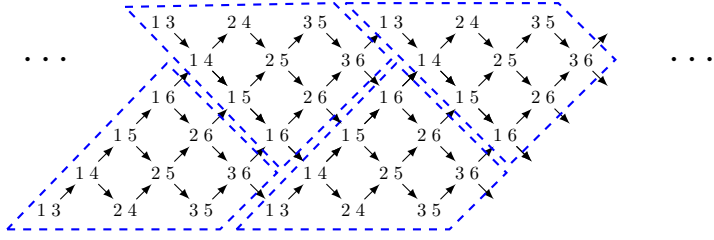

\begin{example}
When
$n=2,t=1$, $\cal A_{2,1}=D^b(\mathbb{K}\vec{A}_{6})/\tau^2[1]$ is  $2$-CY with cluster tilting objects, whose Auslander-Reiten quiver is shown in figure \ref{5}.  By Theorem \ref{d}, we have that the number of torsion pairs in $\cal A_{2,1}: N_{2,1}=20$. We can construct them as follows:
\begin{eqnarray*}
\X_{1}=\{[(0)]\} &\quad\quad& \X^\perp_{1}[-1]=\cal A_{2,1}\nonumber \\
\X_{2}=\{[(13)] \}&\quad\quad& \X^\perp_{2}[-1]=\{[(13)],[(14)],[(16)],[(36)]\}\nonumber \\
\X_{3}=\{[(24)] \}&\quad\quad&\X^\perp_{3}[-1]=\{[(14)],[(15)],[(24)],[(25)]\}\nonumber \\
\X_{4}=\{[(35)] \}&\quad\quad&\X^\perp_{4}[-1]=\{[(25)],[(26)],[(35)],[(36)]\}\nonumber \\
\X_{5}=\{[(14)] \}&\quad\quad& \X^\perp_{5}[-1]=\{[(13)],[(14)],[(24)]\}\nonumber \\
\X_{6}=\{[(25)] \}&\quad\quad& \X^\perp_{6}[-1]=\{[(24)],[(25)],[(35)]\}\nonumber \\
\X_{7}=\{[(36)] \}&\quad\quad& \X^\perp_{7}[-1]=\{[(13)],[(35)],[(36)]\}\nonumber \\
\Y_{1}=\{[(13)],[(14)] \}&\quad\quad& \Y^\perp_{1}[-1]=\{[(13)],[(14)]\}=\Y_{1}\nonumber \\
\Y_{2}=\{[(24)],[(25)] \}&\quad\quad& \Y^\perp_{2}[-1]=\{[(24)],[(25)]\}=\Y_{2}\nonumber \\
\Y_{3}=\{[(35)],[(36)] \}&\quad\quad& \Y^\perp_{3}[-1]=\{[(35)],[(36)]\}=\Y_{3}\nonumber \\
\Y_{4}=\{[(14)],[(24)]\}&\quad\quad& \Y^\perp_{4}[-1]=\{[(14)],[(24)]\}=\Y_{4}\nonumber \\
\Y_{5}=\{[(25)],[(35)] \}&\quad\quad& \Y^\perp_{5}[-1]=\{[(25)],[(35)]\}=\Y_{5}\nonumber \\
\Y_{6}=\{[(13)],[(36)]\}&\quad\quad& \Y^\perp_{6}[-1]=\{[(13)],[(36)]\}=\Y_{6}
\end{eqnarray*}
\end{example}

By Example \ref{c8}, we have known the number of torsion pairs in $\cal A_{2,2}$ is $T_{2+1}=32$. Next, we should count the number of extension closed subcategory $\X$, which consists of indecomposable rigid objects and contains one and only one of indecomposable rigid objects with maximal length.

The rigid indecomposable objects with maximal length in $\cal A_{2,1}$ are $[(1,4)], [(2,5)], [(3,6)]$. The subcategories containing only indecomposable rigid objects and the rigid object $[(1,4)]$ are the followings: $\X_{5}=\{[(14)]\}$, $\X^\perp_{5}[-1]=\{[(13)],[(14)],[(24)]\}$, $\Y_{1}=\{[(13)],[(14)]\}$, and $\Y_{4}=\{[(14)],[(24)]\}$. Then shifting all the subcategories above, we have all the subcategories containing only indecomposable rigid objects and one of the rigid objects $[(i,i+n+1)]$ with maximal length, where $i\in\{1,2,3\}$. The number of such subcategories is $t_{2,1}=3\times 4=12=(2+1)s_{2+2}$.

Then the number of torsion pairs in $\cal A_{2,1}$ is $T_{2+1}-12=20$. It follows that $(\X_i, \X^\perp_{i})$ is a torsion pair in $\cal A_{2,1}$, where $i=1,2,3,\cdots,7$.
Similarly, we can know that $({^\perp\X_{i}},\X_i)$ is a torsion pair in $\cal A_{2,1}$, where $i=1,2,3,\cdots,7$. Moreover, $(\Y_j, \Y^\perp_{j}[-1])=(\Y_j,\Y_j)$ is a cotorsion pair in $\cal A_{2,1}$, where $j=1,2,3,\cdots,6$. Note that $^\perp\Y_j[1]=\Y^\perp_j[-1]$. It follows that $(\Y_j, \Y_{j}[1])$ is a torsion pair in $\cal A_{2,1}$, where $i=1,2,3,\cdots,6$.

\section{Classification of torsion pairs in $\D_{n,t}$}

In this section, we give a classification of  torsion pairs in $\D_{n,t}$ and count the number. Let $u=2t(n+1)$ and $\mathcal C_{D_u}$ be the cluster category of type $D_{u}$. Then $\cal D_{n,t}=D^{b}(\mathbb{K}\vec{D}_{u})/\tau^{n+1}\varphi^{n}$, where $\varphi$ is induced by an automorphism of $D_u$ of order $2$, $n\geq1,t\geq1$. It follows from Lemma 2.9 in \cite{bpr} that there exists a covering functor $\pi\colon\mathcal C_{D_{u}}\rightarrow\cal D_{n,t}$, which is a triangle functor.

Write $F=\tau^{(n+1)}\varphi^n$, then $F: \mathcal C_{D_u}\rightarrow\mathcal C_{D_u}$ is an autoequivalence and $\D_{n,t}$ is the orbit category $\C_{D_u}/F$ (compare Lemma 2.9 in \cite{bpr}). Since $\tau^{-u+1}=[1]$ in $D^{b}(\mathbb{K}\vec{D}_{u})$ by \cite{s}, we have $\tau^{-2t(n+1)}=1=\tau^{2t(n+1)}$ in $\cal C_{D_u}$,  and $\pi$ is a $2t$-covering functor. By Theorem \ref{a4}, we have the following.

\begin{corollary}
There is a bijection between the set of $F$-periodic torsion pairs in $\mathcal C_{D_{u}}$ and the set of torsion pairs in $\cal D_{n,t}$.
\end{corollary}

Let us recall the definition of Ptolemy diagrams of type $D$ and its relation to torsion pairs in the cluster categories of type $D$ based on \cite{hjr3}.
For any $n\geq 1$ we consider a regular $2n$-gon $Q_n$, we label the vertices of $Q_n$ clockwise by $1,2,\ldots 2n$
consecutively. In our arguments below vertices will also be numbered
by some $r\in\mathbb{N}$ which might not be in the range
$1\le r\le 2n$; in this case the numbering of vertices always has
to be taken modulo $2n$.

An {\em arc} is a set $\{i,j\}$ of vertices of
$Q_n$ with $j\not\in \{i-1,i,i+1\}$, i.e.
$i$ and $j$ are different and non-neighboring vertices.
The arcs connecting two opposite vertices $i$ and $i+n$ are called
{\em diameters}. We need two different copies of each of these
diameters and denote them by $\{i,i+n\}_g$ and $\{i,i+n\}_r$, where
$1\le i\le 2n$. The indices should indicate that these diameters are
coloured in the colours green and red, which is a convenient
way to think about and to visualize the diameters.
By a slight abuse of notation, we sometimes omit the indices
and just write $\{i,i+n\}$ for diameters, to avoid cumbersome
definitions or statements.

Any arc in $Q_n$ which is not a diameter is of the form $\{i,j\}$  where
$j\in [i+2,i+n-1]$; here $[i+2,i+n-1]$ stands for the set of vertices
of the $2n$-gon $Q_n$ which are met when going clockwise from
$i+2$ to $i+n-1$ on the boundary of $Q_n$. See figure \ref{11} for an example, for better visibility we draw the red diameters in a wavelike form and the green ones as straight lines.
\begin{figure}
\centering
\begin{tikzpicture}[scale=0.4][auto]
\draw[very thick] (0,0) circle (4) ;

\draw (80:4) edge[thick, color=black!] (-160:4) ;
\draw (-220:4) edge[thick, color=black!] (20:4) ;
\draw (20:4) edge[thick, color=black!] (-100:4) ;
\draw (-160:4) edge[thick, color=black!] (-40:4) ;

\draw (-160:4) edge[thick, dashed, decorate, decoration=snake, color=red!] (20:4);
\draw (-160:4) edge[thick, dashed, color=green!] (20:4);

\draw (-60*1+140:4.7) node {$1$} ;
\draw (-60*2+140:4.8) node {$2$} ;
\draw (-60*3+140:5) node {$3$} ;
\draw (-60*4+140:4.7) node {$4$} ;
\draw (-60*5+140:5.7) node {$5$} ;
\draw (-60*6+140:4.7) node {$6$} ;

\draw (-60*1+140:4) node {$\bullet$} ;
\draw (-60*2+140:4) node {$\bullet$} ;
\draw (-60*3+140:4) node {$\bullet$} ;
\draw (-60*4+140:4) node {$\bullet$} ;
\draw (-60*5+140:4) node {$\bullet$} ;
\draw (-60*6+140:4) node {$\bullet$} ;
\end{tikzpicture}
\caption{$2n=6$}
\label{11}
\end{figure}

Such an arc has a partner arc $\{i+n,j+n\}$ which is
obtained from $\{i,j\}$ by a rotation by 180 degrees. We denote
the pair of arcs $\{\{i,j\},\{i+n,j+n\}\}$ by $\overline{\{i,j\}}$
throughout this section of the paper. The indecomposable objects in $\C_{{D}_n}$ are
in bijection with the union of the set of pairs $\overline{\{i,j\}}$
of non-diameter arcs and the set of diameters $\{i,i+n\}_g$
and $\{i,i+n\}_r$ in two different colours. This bijection extends to subcategories of $\C_{{D}_n}$ closed under direct sums and direct summands and collections of arcs of $2n$-gon $Q_n$.

For the pair of non-diameter
arcs $\overline{\{i,j\}}$ the corresponding
indecomposable object has coordinates $(i,j)$; note that
the coordinates are only determined modulo $n$ so both arcs
$\{i,j\}$ and $\{i+n,j+n\}$ in the pair $\overline{\{i,j\}}$
yield the same coordinate in the Auslander-Reiten quiver of $\C_{D_{n}}$.
The action of $\tau$ on non-diameter arcs is rotation by one vertex, the action of $\tau$ on diameters is rotation by one vertex and changing their colour \cite{hjr3, s}. We note that $\varphi$ acts on diameters by changing their colour and $\varphi=id$  on non-diameter arcs \cite{bpr}.

Back to our consideration, for a coordinate $(i,j)$ corresponding to an indecomposable object of $\C_{D_{u}}$, we denote $[(i,j)]$ by the image under the covering functor $\pi$, then $[(i,j)]$ determines an indecomposable object of $\D_{n,t}$.
For any subcategory $\X$ in $\D_{n,t}$, the preimage under the covering functor $\pi$ corresponds to an $F$-periodic subcategory $\widetilde{\X}=\pi^{-1}(\X)$ in $\C_{D_{u}}$. Moreover, the subcategory $\widetilde{\X}$ corresponds to set of arcs of $2u$-gon $Q_u$ by the discussion above, we still denote the corresponding set of arcs by $\widetilde{\X}$. Fix $F=\tau^{n+1}\varphi^{n}$. We call that the set of arcs $\widetilde{\X}$ is $F$-$periodic$, if the following conditions are satisfied: $(1)$ For each non-diameter arc $(i, j)\in\widetilde{\X}$, also all arcs $(i+(n+1)r, j+(n+1)r)$ \emph{(}modulo $u$\emph{)} for $1\leq r\leq 2t$ are in $\widetilde{\X}$. $(2)$ For each diameter $(i,i+u)\in\widetilde{\X}$,  all diameters $(i+(n+1)r, j+u+(n+1)r)$ \emph{(}modulo $u$\emph{)} with the same colour as $(i,i+u)$, where  $1\leq r\leq 2t$ and $r$ is even, are in $\widetilde{\X}$, and all diameters $(i+(n+1)r, j+u+(n+1)r)$ \emph{(}modulo $u$\emph{)} with the opposite colour as $(i,i+u)$, where $1\leq r\leq 2t$ and $r$ is odd, are in $\widetilde{\X}$. In the rest of this section, we use $(i,j)$ to represent an indecomposable object in $\C_{D_{u}}$ or an arc of $2u$-gon $Q_{u}$ without confusion.

\begin{lemma}
There is a bijection between the following sets:
\begin{enumerate}
  \item[\emph{(1)}]  Subcategory $\X$ of $\cal D_{n,t}$;
  \item[\emph{(2)}] Collection of arcs $ \widetilde{\X}$ of the $2u$-gon $Q_u$ which are $F$-periodic.
\end{enumerate}
\end{lemma}

We recall the definition of Ptolemy diagrams of type $D$ for a $2n$-gon $Q_n$ from \cite{hjr3}.

\begin{definition}\label{c1}
\begin{enumerate}
\item[{(a)}] We say that two non-diameter arcs $\{i,j\}$ and
$\{k,\ell\}$ \emph{cross} precisely if the elements $i,j,k,\ell$ are all
distinct and come in the order $i, k, j, \ell$
when moving around the $2n$-gon $Q_n$ in one direction or the other
(i.e. counterclockwise or clockwise). In particular, the two
arcs in $\overline{\{i,j\}}$ do not cross.

Similarly, in the case $j=i+n$, the above condition defines
when a diameter $\{i,i+n\}_g$ (or $\{i,i+n\}_r$) crosses
the non-diameter arc $\{k,\ell\}$.
\item[{(b)}] We say that two pairs $\overline{\{i,j\}}$ and
$\overline{\{k,\ell\}}$ of non-diameter arcs {\em cross} if
there exist two arcs in these two pairs which cross in the sense of
part (a). (Note that then necessarily the other two rotated arcs also
cross.)

Similarly, the diameter $\{i,i+n\}_g$ (or $\{i,i+n\}_r$) crosses
the pair $\overline{\{k,\ell\}}$ of non-diameter arcs if it
crosses one of the arcs in $\overline{\{k,\ell\}}$.
(Note that it then necessarily crosses both arcs in $\overline{\{k,\ell\}}$.)
\item[{(c)}] Two diameters $\{i,i+n\}_g$ and $\{j,j+n\}_r$
of different colour {\em cross} if $j\not\in \{i,i+n\}$, i.e.
if they have different endpoints. But $\{i,i+n\}_g$ and
$\{i,i+n\}_r$ do not cross.
Moreover, any diameters of the same colour do not cross.
\end{enumerate}
\end{definition}

\begin{definition}\label{c2}
  Let $\X$ be a collection of arcs of the $2n$-gon $Q_n$, $n> 1$,
  which is invariant
  under rotation of $180$ degrees.  Then $\X$ is called a
  \emph{Ptolemy diagram of type~$D$ }if it satisfies the following
  conditions.  Let $\alpha = \{i,j \}$ and $\beta = \{k,\ell \}$
  be crossing arcs in $\X$ (in the sense of
  Definition \ref{c1}).
  \begin{enumerate}
  \item[{(Pt1)}] If $\alpha$ and $\beta$ are not diameters, then
  those of $\{
    i,k\}$, $\{i,\ell \}$, $\{
    j,k \}$, $\{ j,\ell \}$ which are
    arcs in $P$ are also in $\X$.
    In particular, if two of the vertices $i,j,k,\ell$ are
    opposite vertices (i.e. one of $k$ and $\ell$ is equal to
    $i+n$ or $j+n$), then
    both the green and the red diameter connecting them are also
    in $\X$.
  \item[{(Pt2)}] If both $\alpha$ and $\beta$ are diameters
  (necessarily of different colour by Definition \ref{c1}\,(c))
  then those of $\{ i,k\}$, $\{
    i,k+n \}$, $\{ i+n,k \}$, $\{
    i+n,k+n \}$ which are arcs of $P$ are also in $\X$.
\item[{(Pt3)}] If $\alpha$ is a diameter while $\beta$ is not a diameter,
  then those of $\{
    i,k\}$, $\{ i,\ell \}$, $\{
    j,k \}$, $\{ j,\ell \}$ which are
    arcs and do not cross the arc $\{k+n,\ell +n \}$ are also in
    $\X$.
    Additionally, the diameters $\{ k, k+n \}$ and
    $\{ \ell, \ell +n \}$ of the same colour as $\alpha$ are
    also in $\X$.
  \end{enumerate}

\end{definition}
These conditions are illustrated in figure \ref{10}.
\begin{figure}
  \centering
  \begin{enumerate}
  \item[{(Pt1)}]
  The first Ptolemy condition in type
    $D$:
$$
\begin{tikzpicture}[scale=0.4]
\draw[very thick] (0,0) circle (4) ;
\draw (90:4) edge[thick, dashed, color=black!] (45:4) ;
\draw (45:4) edge[thick, dashed, color=black!] (-180:4) ;
\draw (-180:4) edge[thick, dashed, color=black!] (-225:4) ;
\draw (-225:4) edge[thick, dashed, color=black!] (90:4) ;
\draw (90:4) edge[thick, color=black!] (-180:4) ;
\draw (45:4) edge[thick, color=black!] (-225:4) ;

\draw (0:4) edge[thick, dashed, color=black!] (-45:4) ;
\draw (-45:4) edge[thick, dashed, color=black!] (-90:4) ;
\draw (-90:4) edge[thick, dashed, color=black!] (-135:4) ;
\draw (-135:4) edge[thick, dashed, color=black!] (0:4) ;
\draw (0:4) edge[thick, color=black!] (-90:4) ;
\draw (-45:4) edge[thick, color=black!] (-135:4) ;

\draw (-45*1+135:4.7) node {$\ell$} ;
\draw (-45*2+135:4.7) node {$j$} ;
\draw (-45*3+135:5.3) node {$k+n$} ;
\draw (-45*4+135:5) node {$i+n$} ;
\draw (-45*5+135:4.7) node {$\ell+n$} ;
\draw (-45*6+135:5) node {$j+n$} ;
\draw (-45*7+135:4.7) node {$k$} ;
\draw (-45*8+135:4.7) node {$i$} ;

\draw (-45*1+135:4) node {$\bullet$} ;
\draw (-45*2+135:4) node {$\bullet$} ;
\draw (-45*3+135:4) node {$\bullet$} ;
\draw (-45*4+135:4) node {$\bullet$} ;
\draw (-45*5+135:4) node {$\bullet$} ;
\draw (-45*6+135:4) node {$\bullet$} ;
\draw (-45*7+135:4) node {$\bullet$} ;
\draw (-45*8+135:4) node {$\bullet$} ;
\end{tikzpicture}
\begin{tikzpicture}[scale=0.4][auto]
\draw[very thick] (0,0) circle (4) ;
\foreach \a in {1,2,...,5}{
\draw (-60*\a+140:4) edge[thick, dashed, color=black!] (-60*\a+80:4) ;
}
\draw (80:4) edge[thick, dashed, color=black!] (-220:4) ;

\draw (80:4) edge[thick, color=black!] (-160:4) ;
\draw (-220:4) edge[thick, color=black!] (20:4) ;
\draw (20:4) edge[thick, color=black!] (-100:4) ;
\draw (-160:4) edge[thick, color=black!] (-40:4) ;

\draw (-160:4) edge[thick, dashed, decorate, decoration=snake, color=red!] (20:4);
\draw (-160:4) edge[thick, dashed, color=green!] (20:4);

\draw (-60*1+140:4.7) node {$\ell$} ;
\draw (-60*2+140:4.8) node {$\;\;\;\;\;\;\;\;j=k+n$} ;
\draw (-60*3+140:5) node {$i+n$} ;
\draw (-60*4+140:4.7) node {$\ell+n$} ;
\draw (-60*5+140:5.7) node {$k=j+n$} ;
\draw (-60*6+140:4.7) node {$i$} ;

\draw (-60*1+140:4) node {$\bullet$} ;
\draw (-60*2+140:4) node {$\bullet$} ;
\draw (-60*3+140:4) node {$\bullet$} ;
\draw (-60*4+140:4) node {$\bullet$} ;
\draw (-60*5+140:4) node {$\bullet$} ;
\draw (-60*6+140:4) node {$\bullet$} ;
\end{tikzpicture}
$$

  \item[{(Pt2)}]
The second Ptolemy condition in type
  $D$:
  $$
\begin{tikzpicture}[scale=0.4]
\draw[very thick] (0,0) circle (4) ;

\draw [thick, dashed, color=black!] (130:4) -- (70:4) ;
\draw [thick, dashed, color=black!] (70:4) -- (-50:4) ;
\draw [thick, dashed, color=black!] (-50:4) -- (-110:4) ;
\draw [thick, dashed, color=black!] (-110:4) -- (130:4) ;
\draw [thick, color=green!] (70:4) -- (-110:4);
\draw [thick, decorate, decoration=snake, color=red!] (130:4) -- (-50:4);

\draw (130:4.7) node {$i$} ;
\draw (70:4.7) node {$k+n$} ;
\draw (-50:5) node {$i+n$} ;
\draw (-110:4.7) node {$k$} ;

\draw (130:4) node {$\bullet$} ;
\draw (70:4) node {$\bullet$} ;
\draw (-50:4) node {$\bullet$} ;
\draw (-110:4) node {$\bullet$} ;
\end{tikzpicture}
$$

\item[{(Pt3)}]
The third Ptolemy condition in type
  $D$:
\[
\begin{tikzpicture}[scale=0.4]
\draw[very thick] (0,0) circle (4);
\draw [thick, dashed, color=black!] (60:4) -- (120:4);
\draw [thick, dashed, color=black!] (120:4) -- (160:4);
\draw [thick, dashed, color=black!] (-120:4) -- (-60:4);
\draw [thick, dashed, color=black!] (-60:4) -- (-20:4);

\draw [very thick, color=black!] (160:4) -- (60:4);
\draw [very thick, color=black!] (-120:4) -- (-20:4);

\draw [thick, color=green!] (120:4) -- (-60:4);
\draw [thick, dashed, color=green!] (60:4) -- (-120:4);
\draw [thick, dashed, color=green!] (160:4) -- (-20:4);

\draw (60:4.7) node {$\ell$} ;
\draw (120:4.7) node {$i$} ;
\draw (160:4.7) node {$k$} ;
\draw (-120:4.7) node {$\ell+n$} ;
\draw (-60:4.8) node {$i+n$} ;
\draw (-20:5) node {$k+n$} ;

\draw (60:4) node {$\bullet$} ;
\draw (120:4) node {$\bullet$} ;
\draw (160:4) node {$\bullet$} ;
\draw (-120:4) node {$\bullet$} ;
\draw (-60:4) node {$\bullet$} ;
\draw (-20:4) node {$\bullet$} ;
\end{tikzpicture}
\]
\end{enumerate}

\caption{The Ptolemy conditions in type~$D$.}
\label{10}
\end{figure}

Now we define the $F$-$periodic$ $Ptolemy$ $diagram$ of type $D$ for $2u$-gon $Q_u$ which we will use to give a classification of torsion pairs in $\cal D_{n,t}$.

\begin{definition} Let $\widetilde{\X}$ be a collection of arcs in $2u$-gon $Q_{u}$, and $F=\tau^{n+1}\varphi^{n}$. $\widetilde{\X}$ is called an $F$-$periodic$ $Ptolemy$ $diagram$ of type $D$ if $\widetilde{\X}$ is a Ptolemy diagram of type $D$ and is $F$-periodic.
\end{definition}

\begin{definition}\label{c3}
For a coordinate $(i,j)$ (modulo $2u$) with $j>i$ corresponding to an indecomposable object in $\mathcal C_{D_{u}}$, we call $j-i-1$ the $level$ of the vertex, and $j-i$ the $length$ of the vertex. For an indecomposable object $[(i,j)]$ in $\D_{n,t}$, we also call $j-i-1$ the $level$ of the vertex. The $length$ of the vertex $[(i,j)]$ is defined by $j-i$.
\end{definition}

\begin{lemma}\label{c4}
Let $(\X,\X^{\perp})$ be a torsion pair in $\cal D_{n,t}$ and $\widetilde{\X}$ be the corresponding
$F$-periodic collection of arcs of the $2u$-gon $Q_u$, where $u=2t(n+1)$. If $t>1$, then precisely one
of the following situation occurs:
\begin{enumerate}
\item[{\emph{(i)}}] The level of all  indecomposable objects of $\X\leq n$;

\item[\emph{(ii)}] The level of all  indecomposable objects of $\X^{\perp}\leq n$.
\end{enumerate}
\end{lemma}

\proof Recall that the elements with level $\leq n$ in $\cal D_{n,t}$ are exactly all the rigid indecomposable
objects \cite[Lemma 2.9]{bpr}.
\begin{enumerate}

\item[(I)] If the level of all the indecomposable objects of $\X\leq n$,  we claim that $\X^{\bot}$ must contain an element (not diameters) with level$>n$. Indeed, since $\X$ contains only (finitely many) indecomposable rigid objects, we pick an arc from $\X$ with maximal length. We can suppose that its coordinate is $[(1,\ell)]$ with $3\leq \ell\leq n+2$, since $[(1,\ell)]$ corresponds to a rigid indecomposable object. If we can show $[(1,\ell+n+1)]\in\X^{\bot}[-1]$, then $\X^{\bot}$  contains an element with level$>n$. Firstly, since $t>1$, $\ell+n+1\leq n+2+n+1=2n+3<u+1$, that is $[(1,\ell+n+1)]$ represents a non-diameter arc, and the level of $[(1,\ell+n+1)]$ is $\ell+n-1\geq 3+n-1=n+2$. Secondly, since $(1,\ell)$ is in $\widetilde{\X}$ with maximal length, there is no arc $(i,j)\in\widetilde{\X}$ with $1<i<\ell<j$. Otherwise $(1,\ell)$ has to cross $(i,j)$, but the Ptolemy condition yields an arc longer than $(1,\ell)$, a contradiction. This means $(1,\ell+n+1)$ cannot cross any arc from $\widetilde{\X}$. Similarly, we can prove $(n+2,\ell+2n+2),(2n+3,\ell+3n+3),\ldots \in\widetilde{\X}$, i.e., $[(1,\ell+n+1)]\in\X^{\bot}[-1]$,  i.e., $(1,\ell+n+1)\in\X^{\bot}[-1]$.

\item[(II)] Suppose $\cal X$ contains an element with level$>n$. Note that  all indecomposable
objects of $\cal D_{n,t}$ with level$>n$ may diameters  or non-diameters, we claim that
$\X^{\perp}$ contains only elements with level$\leq n$.
\item[1)] If $\X$ contains a non-diameter arc with level$>n$, without losing generality, we suppose $[(1,\ell)]$ is such an element. Note that we can choose $n+3\leq\ell\leq u$ \cite[Fig. 6]{bpr}. We choose $\ell$ for different intervals $[n+3,2n+3]$, $[2n+3,3n+4]$,$[3n+4,4n+5],\ldots, [(2t-1)(n+1)+1,2t(n+1)+1]$.

If $n+3\leq\ell\leq 2n+3$, then the corresponding arcs in $\widetilde{\X}$ are shown in figure \ref{6}.
\begin{figure}
\centering
\begin{tikzpicture}[scale=0.4]
\draw[very thick] (0,0) circle (4) ;
\foreach \a in {1,3,...,11}
\draw (-30*\a+120:4) edge[very thick, color=black!40, out={-40-30*\a}, in={190-30*\a}] (-30*\a+30:4) ;

\foreach \x in {1,2,...,12} {
\draw (-30*\x+120:4) node {$\bullet$} ;
} ;
\draw (-30*1+120:4.7) node {$1$} ;
\draw (-30*2+120:4.7) node {} ;
\draw (-30*3+120:5) node {\quad$n+2$} ;
\draw (-30*4+120:4.7) node {$\ell$} ;
\draw (-30*5+120:5.1) node {\quad$2n+3$} ;
\draw (-30*6+120:5) node {\qquad$\ell+n+1$} ;
\draw (-30*7+120:4.7) node {$3n+4$} ;
\draw (-30*8+120:5) node {$\ell+2n+2$\;\;\;\;\;\;\;\;\;\;} ;
\draw (-30*9+120:5.3) node {$4n+5$\;\;\;} ;
\draw (-30*10+120:5.9) node {$\ell+3n+3$\;\;} ;
\draw (-30*11+120:4.7) node {} ;
\draw (-30*12+120:4.7) node {} ;
\end{tikzpicture}
\caption{A case when $t=3$}
\label{6}
\end{figure}

If $2n+3\leq\ell\leq 3n+4$, then the corresponding arcs in $\widetilde{\X}$ are shown in figure \ref{7}.
\begin{figure}
\centering
\begin{tikzpicture}[scale=0.4]
\draw[very thick] (0,0) circle (4) ;
\foreach \a in {1,3,8,10}
\draw (-25.7*\a+115.7:4) edge[very thick, color=black!40, out={-40-25.7*\a}, in={170-25.7*\a}] (-25.7*\a+12.9:4) ;
\foreach \a in {4,6}
\draw (-25.7*\a+115.7:4) edge[very thick, color=black!40, out={-47.3-25.7*\a}, in={155.7-25.7*\a}] (-25.7*\a-12.8:4) ;
\foreach \a in {13}
\draw (-25.7*\a+115.7:4) edge[very thick, color=black!40, out={-37.3-25.7*\a}, in={195.7-25.7*\a}] (-25.7*\a+38.6:4) ;

\foreach \x in {1,2,...,14} {
\draw (-25.7*\x+115.7:4) node {$\bullet$} ;
} ;
\draw (-25.7*1+115.7:4.7) node {$1$} ;
\draw (-25.7*2+115.7:4.7) node {} ;
\draw (-25.7*3+115.7:5) node {\quad$n+2$} ;
\draw (-25.7*4+115.7:5.3) node {\quad$2n+3$} ;
\draw (-25.7*5+115.7:4.7) node {$\ell$} ;
\draw (-25.7*6+115.7:5) node {\quad$3n+4$} ;
\draw (-25.7*7+115.7:4.8) node {\;\;\qquad$\ell+n+1$} ;
\draw (-25.7*8+115.7:5) node {$4n+5$} ;
\draw (-25.7*9+115.7:4.8) node {$\ell+2n+2$\;\;\;\;\;\;\;\;\;\;\;\;} ;
\draw (-25.7*10+115.7:5) node {$5n+6$\;\;\;\;} ;
\draw (-25.7*11+115.7:5.8) node {$\ell+3n+3$\;\;\;} ;
\draw (-25.7*12+115.7:5.8) node {$\ell+4n+4$\;\;\;\;} ;
\draw (-25.7*13+115.7:4.7) node {} ;
\draw (-25.7*14+115.7:4.8) node {$\ell+5n+6$\;\;\;\;\;\;\;\;\;\;\;\;\;\;\;} ;
\end{tikzpicture}
\caption{A case when $t=3$}
\label{7}
\end{figure}

The other cases are similar. This shows the level of indecomposable objects in $[(1,\ell)]^{\perp}[-1]\leq n$. Therefore the level of indecomposable objects in $\X^{\perp}\leq n$, since $\X^{\perp}[-1]\subseteq[(1,\ell)]^{\perp}[-1]$.

\item[2)] If $\X$ contains a diameter $[(1,u+1)]$ (green or red) in $\X$, without losing generality,
we suppose $[(1,u+1)]_g\in\X$. Then $(1,u+1)_g,(n+2,u+n+2)_r,\cdots,\in\widetilde{\X}$
and the first two objects cross. Then $\mathrm{(Pt2)}$ implies that
$(n+2,u+1)\in\widetilde{\X}$ which is not a diameter and then the non-diameter  $[(n+2,u+1)]\in\X$,  its
level is $u+1-(n+2)-1=u-(n+2)=(2t-1)(n+1)-1\geq 3n+2>n$, since $t>1$. Thus $\X^{\perp}$ contains only indecomposable rigid objects by  case $1)$.
\end{enumerate}
As a consequence, for a torsion pair $(\X, \X^{\bot})$  in $\cal D_{n,t}$, precisely one of $(1)$ and $(2)$ occurs. We complete the proof.\qed

\begin{definition}
Let $(i,j)$ be a non-diameter arc of $2u$-gon $Q_{u}$. The $wing$ $W(i,j)$ of $(i,j)$ consists of all arcs $(r,s)$ of the $2u$-gon such that $i\leq r\leq s\leq j$, that is all arcs which are overarched by $(i,j)$.  $[(i,j)]$ represents a vertex in the AR-quiver of $\cal D_{n,t}$, the corresponding wing is denoted by $W[(i,j)]$.
\end{definition}

 Combining Lemma \ref{c4}, we have the classification of torsion pairs in $\cal D_{n,t}, t>1$, whose proof is the same as Theorem \ref{b} (compare \cite{hjr2}).
\begin{theorem}
There are bijections between the following sets for $t>1$:
\begin{enumerate}
  \item[\emph{(1)}] Torsion pairs $(\X, \X^{\bot})$ in $\cal D_{n,t}$ such that the level of all the indecomposable objects in $\X\leq n$;
  \item[\emph{(2)}] $F$-periodic Ptolemy diagrams $\widetilde{\X}$ of type $D$ of $2u$-gon $Q_{u}$ such that all arcs in $\widetilde{\X}$ have length at most $n+1$;
  \item[\emph{(3)}] Collections $\left\{([(i_{1},j_{1})], [W_{1}]), \ldots, ([(i_{r},j_{r})], [W_{r}])\right\}$ of pairs consisting of vertices $[(i_{\ell},j_{\ell})]$ of level$\leq n$ in the AR-quiver of $D_{n,t}$ and subset $[W_{\ell}]\subset W[(i_{\ell},j_{\ell})]$ of their wings such that for any different $k, \ell\in \left\{1,2,\ldots,r\right\}$, we have
  $$ W[(i_{k},j_{k})][1]\cap W[(i_{\ell},j_{\ell})]=\emptyset,$$
  and the $F$-periodic collection $W_{\ell}$ corresponding to $[W_{\ell}]$ is a Ptolemy diagram of type $D$.

\end{enumerate}
\label{f}
\end{theorem}
Next, we describe torsion pairs in $\D_{n,1}$. Recall that two diameters are called $paired$ if they connect the same two vertices (and thus of different colour).
\begin{theorem} For a torsion pair $(\X, \Y)$ in $\D_{n,1}$, precisely one of the following situations occurs:
\begin{itemize}
  \item [1.] $\X$ (resp. $\Y$) contains one paired diameters and $\Y$ (resp. $\X$) contains only indecomposable rigid objects.
  \item [2.] Both $\X$ and $\Y$ contain only one non-paired diameter and some indecomposable rigid objects.
\end{itemize}
\end{theorem}

\proof  Let $(\X,\Y)$ be a torsion pair in $\cal D_{n,1}$ and $\widetilde{\X}$, $\widetilde{\Y}$ be the corresponding
$F$-periodic collections of arcs of the $2u$-gon $Q_{u}$, where $u=2(n+1)$. We note that $(\X,\Y)$ is a torsion pair if and only if so is $(\Y, \X[2])$, since $\cal C$ is $2$-CY.

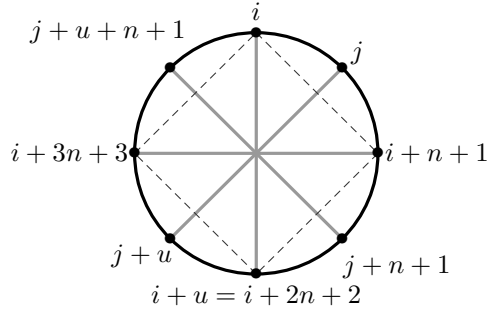
\begin{figure}
\centering
\begin{tikzpicture}[scale=0.4]
\draw[very thick] (0,0) circle (4) ;
\foreach \a in {1,2,3,4}\draw[very thick, color=black!40] (-45*\a+90:4) -- (-45*\a-90:4) ;
\foreach \a in {1,3,5,7}\draw[densely dashed, color=black!80] (-45*\a+135:4) -- (-45*\a+45:4) ;

\foreach \x in {1,2,...,8} {
\draw (-45*\x+135:4) node {$\bullet$} ;
} ;
\draw (-45*1+135:4.7) node {$i$} ;
\draw (-45*2+135:4.7) node {$j$} ;
\draw (-45*3+135:5.7) node {\;\;$i+n+1$} ;
\draw (-45*4+135:5.3) node {\qquad$j+n+1$} ;
\draw (-45*5+135:4.7) node {$i+u=i+2n+2$} ;
\draw (-45*6+135:4.8) node {$j+u$\;\;\;} ;
\draw (-45*7+135:5.9) node {$i+3n+3$\;\;} ;
\draw (-45*8+135:5.7) node {$j+u+n+1$\;\;\;\;\;\;\;} ;
\end{tikzpicture}
\caption{The case when $t=1$}
\label{8}
\end{figure}

\begin{itemize}
\item [(1)] If $\widetilde{\X}$ contains a non-diameter arc with length longer than $n+1$, then $\Y$ contains only indecomposable rigid objects. The proof is the same as the part (II) of the proof of Lemma \ref{c4}.
\item[(2)] If $\widetilde{\X}$ contains a paired diameters, then the indecomposable objects in $\Y$ are all rigid. Suppose
$[(1,u+1)]\in\widetilde{\X}$ (green and red) is one paired diameters. Because $\widetilde{\X}$ is $F$-periodic, $(1,u+1)_{r,g}\in\widetilde{\X}$, $(n+2,3n+4)_{r,g}\in\widetilde{\X}$ and they cross. So $(1,n+2)$ and $(n+2,2n+3)$ are in $\widetilde{\X}$ by $\mathrm{(Pt2)}$ of Definition \ref{c2}. Thus the maximal length of arcs that do not cross  any arc in $\widetilde{\X}$ is $n+1$ and any diameter crosses $(1,u+1)$ in $\widetilde{\X}$, that is $\X^{\perp}[-1]$ contains only indecomposable rigid objects, so does $\Y$.
\item[(3)] If $\X$ contains only indecomposable rigid objects, then $\widetilde{\Y}$ will contain a paired diameters. Suppose the arc in $\widetilde{\X}$ with maximal length is $(1,\ell)$ (up to shifting), then $3\leq \ell\leq n+2$ since $(1,\ell)$ corresponds to a rigid object. We claim the paired diameters $(2,u+2)$ are in $\widetilde{\Y}$.
    Otherwise, there is an arc in $\widetilde{\X}$ crosses the paired diameters $(1,u+1)$, then there exists an arc with length longer than $(1,\ell)$ similarly as the proof of Lemma \ref{c4}. So the paired diameters $(1,u+1)$ are in $\Y[-1]$, that is, $(2,u+2)$ (red and green) are in $\Y$.
\item[(4)]If $\X$ contains diameters, but no paired diameter, then $\widetilde{\X}$ contains no paired diameters. We first claim that $\X$ contains only one diameter (red or green). In fact, if $\X$ contains $2$ non-paired diameters $[(i,i+u)],[(j,j+u)]$ (red or green) $(j\neq i,j\neq i+n+1)$, then $\widetilde{\X}$ will contain $4$ non-paired diameters with different colours, see figure \ref{8}.

    Note that $(i,i+u)$ and $(i+n+1,i+3n+3)$ are two different coloured diameters and they cross, then $(i,i+n+1)\in\widetilde{\X}$ by (Pt2). Moreover, the arc $(i,i+n+1)$ crosses  $(j,j+u)$, then $(i,i+u)$ and $(i+n+1,i+3n+3)$ with the same colour as $(j,j+u)$ are in $\widetilde{\X}$, but $(i,i+u)$ and $(i+n+1,i+3n+3)$ are different colours, this implies
that $[(i,i+u)]$ are paired diameters in $\X$, a contradiction.

Moreover,  if $\X$ contains only one diameter, without losing generality, we assume its coordinate is $[(1,1+u)]$. Then $(1,1+u)$ and $(n+2,u+n+1)$ are in $\widetilde{\X}$ and they are different colours, so they cross. The Ptolemy condition implies that $(1,n+2)$ and $(n+2,2n+3)$ are in $\widetilde{\X}$, so the maximal length of arcs that do not cross any arc in $\widetilde{\X}$ is $n+1$, and the diameter $[(1,1+u)]$ with different colour as it in $\X$ does not cross any arc in $\widetilde{\X}$ either. That means that $\X^{\perp}[-1]$ contains only indecomposable rigid objects and one diameter, so does $\Y=\X^{\perp}$.

    \end{itemize}

As a consequence, if $\X$ (resp. $\Y$) contains only indecomposable rigid objects, then (3) ensures $\Y$ (resp. $\X$) contains a paired diameters, so case 1 occurs. Suppose $\X$ (resp. $\Y$) contains a non-rigid object. If the non-rigid object is non-diameter, then (1) ensures $\Y$ (resp. $\X$) contains only indecomposable rigid objects, and case 1 occurs. If the non-rigid object is a diameter, then (2) implies that case 1 occurs if the diameters are paired, and  (4) implies that case 2 occurs if the diameter is non-paired. Obviously, case 1 and case 2 cannot occur simultaneously.
\qed

\begin{theorem} Let $D_{n,t}$ be the number of torsion pairs in $\cal D_{n,t}$.
 \begin{itemize}
   \item[1] If $t>1$, then $D_{n,t}=T_{n+1},$  the number of torsion pairs in $\cal A_{n,t}$.
   \item[2] If $t=1$, then $D_{n,t}=T_{n+1}+2t_{n,1}=\sum\limits_{\ell\geq 0}2^{\ell+1}\binom{n+\ell}{\ell}\bigg[\binom{2n+1}{n-2\ell}+
       \binom{2n}{n-2\ell}\bigg]$
 \end{itemize}
\label{h2}
\end{theorem}

\proof For $t>1$, we only consider the subcategories $\X$ of $\cal D_{n,t}$ with the level of indecomposable objects in
$\X\leq n$. Because such subcategories $\X$ cannot contain any diameter and $\mathrm{(Pt1)}$ in Definition \ref{c2} coincides with the Ptolemy condition in type $A$, the Ptolemy diagram of type $D$ is the same as Ptolemy diagram of type $A$. Moreover, because $\X$ contains only indecomposable rigid objects, the $F$-periodic of the corresponding collection of non-diameter arcs are $(n+1)$-periodic.  By  Theorem 4.9 and Theorem \ref{b}, we know that
the number of torsion pairs in $\cal D_{n,t}$ equals to
the number of torsion pairs in $\cal A_{n,t}$ and equals to
the number of torsion pairs in the cluster tube of rank $n+1$.

For $t=1$, The torsion pairs in $\cal D_{n,1}$ divide into two subclasses:  one is the torsion pairs $(\X,\Y)$ in $\cal D_{n,1}$ with $\X$ or $\Y$ (not both) containing a paired diameters, counting the number of this case reduces to counting of the possible halves $\X$ or $\Y$ of a torsion pair, whose all indecomposable objects are rigid by Theorem 4.10. The number of such torsion pairs in $\cal D_{n,1}$ equals to the number of torsion pairs in $\A_{n,t}$ with $t>1$. Another one is the torsion pairs $(\X,\Y)$ in $\cal D_{n,1}$ with both $\X$ and $\Y$ containing one non-paired diameter and some indecomposable rigid objects. The number of such torsion pairs is $2t_{n,1}$.\qed

\begin{example}
When
$n=1,t=2$, $\cal D_{1,2}=D^b(\mathbb{K}D_{8})/\tau^2\varphi$ is  $2$-Calabi-Yau with maximal rigid objects, whose Auslander-Reiten quiver is shown in  figure \ref{13}. By Theorem \ref{h2}, the number of torsion pairs in $\cal D_{1,2}$ is $T_{1+1}=6$. We list the torsion pairs according to Theorem \ref{f}:
\begin{figure}
\centering
\begin{tikzpicture}[scale=0.75,
fl/.style={->,shorten <=6pt, shorten >=6pt,>=latex}]
\foreach \x in {1,2} {
 \pgfmathparse{8}\let\z\pgfmathresult ;
  \foreach \y in {2,3,...,7} {
   \newcount\u ;
   \pgfmathsetcount{\u}{\x+\y} ;
   \draw (2*\x+\y-4,\y-2) node[scale=0.7] {\x$\;$\the\u} ;
  } ;
   \draw (2*1+8-4,8-2) node[scale=0.7] {$1\;9^+$} ;
   \draw (2*2+8-4,8-2) node[scale=0.7] {$2\;10^+$} ;
  \foreach \y in {2,3,...,7} {
   \newcount\u ;
   \pgfmathsetcount{\u}{\x+\y} ;
   \draw[fl] (\y-4+2*\x,\y-2) -- (\y-3+2*\x,\y-1) ;
   \draw[fl] (\y-3+2*\x,\y-1) -- (\y-2+2*\x,\y-2) ;
  } ;
  \draw[fl] (5,5) -- (6,5) ;
  \draw (6,5) node[scale=0.7] {$1\;9^-$} ;
  \draw[fl] (6,5) -- (7,5) ;
  \draw[fl] (7,5) -- (7.9,5) ;
  \draw (8,5) node[scale=0.7] {$2\;10^-$} ;
  \draw[fl] (8,5) -- (9,5) ;
} ;
\draw[thick, dashed, blue] (-1.1,-0.5) -- ++(4,0) -- ++(6.8,6.8) -- ++(-4,0) -- cycle ;

\begin{scope}[xshift=8cm]
\foreach \x in {1,2} {
 \pgfmathparse{8}\let\z\pgfmathresult ;
  \foreach \y in {2,3,...,7} {
   \newcount\u ;
   \pgfmathsetcount{\u}{\x+\y} ;
   \draw (2*\x+\y-8,\y-2) node[scale=0.7] {\x$\;$\the\u} ;
  } ;
   \draw (2*1+8-8,8-2) node[scale=0.7] {$1\;9^+$} ;
   \draw (2*2+8-8,8-2) node[scale=0.7] {$2\;10^+$} ;
  \foreach \y in {2,3,...,7} {
   \newcount\u ;
   \pgfmathsetcount{\u}{\x+\y} ;
   \draw[fl] (\y-8+2*\x,\y-2) -- (\y-7+2*\x,\y-1) ;
   \draw[fl] (\y-7+2*\x,\y-1) -- (\y-6+2*\x,\y-2) ;
  } ;
  \draw[fl] (1,5) -- (2,5) ;
  \draw (2,5) node[scale=0.7] {$1\;9^-$} ;
  \draw[fl] (2,5) -- (3,5) ;
  \draw[fl] (3,5) -- (3.9,5) ;
  \draw (4,5) node[scale=0.7] {$2\;10^-$} ;
  \draw[fl] (4,5) -- (5,5) ;
} ;
\end{scope}
\begin{scope}[xshift=240]
\draw[thick, dashed, blue] (-5.3,-0.5) -- ++(3.7,0) -- ++(6.8,6.8) -- ++(-3.7,0) -- cycle ;
\end{scope}
\draw (0,3) node[scale=1.5] {$\cdots$} ;
\draw (13,3) node[scale=1.5] {$\cdots$} ;
\end{tikzpicture}
\caption{The AR-quiver of $\cal D_{1,2}$}
\label{13}
\end{figure}
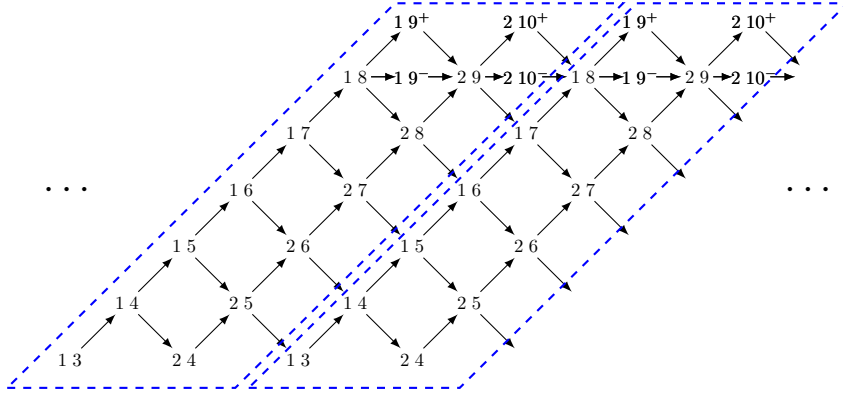
\begin{eqnarray*}
\X_{1}=\{[(0)] \}&\quad\quad& \X^\perp_{1}[-1]=\cal D_{1,2}\nonumber \\
\X_{2}=\{[(13)] \}&\quad\quad& \X^\perp_{2}[-1]=\{[(13)],[(15)],[(17)],[( 19^{+})],[(19^{-})]\}\nonumber \\
\X_{3}=\{[(24)] \}&\quad\quad&\X^\perp_{3}[-1]=\{[(24)],[(26)],[(28)],[(2 10^{+})],[(2 10^{-})]\}
\end{eqnarray*}
\end{example}
We only need to consider the subcategories containing only indecomposable rigid objects in $\cal D_{1,2}$. By Theorem \ref{f}, all the torsion pairs $(\X_{i}, \X^\perp_{i})$ and $( ^\perp\X_{i}, \X_{i})$ are listed above for $i=1,2,3$.

\begin{example}
When
$n=1,t=1$, $\cal D_{1,1}=D^b(\mathbb{K}D_{4})/\tau^2\varphi$ is  $2$-Calabi-Yau with maximal rigid objects, whose Auslander-Reiten quiver is shown in figure \ref{12}. By Theorem \ref{h2}, the number of torsion pairs in $\cal D_{1,1}$ is $T_{1+1}+2t_{1,1}=10$. We list the torsion pairs according to Theorem 4.10:

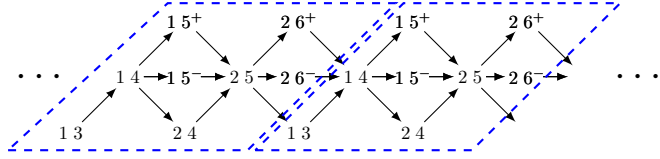
\begin{figure}
\centering
\begin{tikzpicture}[scale=0.75,
fl/.style={->,shorten <=6pt, shorten >=6pt,>=latex}]
\foreach \x in {1,2} {
 \pgfmathparse{4}\let\z\pgfmathresult ;
  \foreach \y in {2,3} {
   \newcount\u ;
   \pgfmathsetcount{\u}{\x+\y} ;
   \draw (2*\x+\y-4,\y-2) node[scale=0.7] {\x$\;$\the\u} ;
  } ;
   \draw (2*1+4-4,4-2) node[scale=0.7] {$1\;5^+$} ;
   \draw (2*2+4-4,4-2) node[scale=0.7] {$2\;6^+$} ;
  \foreach \y in {2,3} {
   \newcount\u ;
   \pgfmathsetcount{\u}{\x+\y} ;
   \draw[fl] (\y-4+2*\x,\y-2) -- (\y-3+2*\x,\y-1) ;
   \draw[fl] (\y-3+2*\x,\y-1) -- (\y-2+2*\x,\y-2) ;
  } ;
  \draw[fl] (1,1) -- (2,1) ;
  \draw (2,1) node[scale=0.7] {$1\;5^-$} ;
  \draw[fl] (2,1) -- (3,1) ;
  \draw[fl] (3,1) -- (3.9,1) ;
  \draw (4,1) node[scale=0.7] {$2\;6^-$} ;
  \draw[fl] (4,1) -- (5,1) ;
} ;
\draw[thick, dashed, blue] (-1.1,-0.3) -- ++(4.2,0) -- ++(2.6,2.6) -- ++(-4,0) -- cycle ;

\begin{scope}[xshift=8cm]
\foreach \x in {1,2} {
 \pgfmathparse{4}\let\z\pgfmathresult ;
  \foreach \y in {2,3} {
   \newcount\u ;
   \pgfmathsetcount{\u}{\x+\y} ;
   \draw (2*\x+\y-8,\y-2) node[scale=0.7] {\x$\;$\the\u} ;
  } ;
   \draw (2*1+4-8,4-2) node[scale=0.7] {$1\;5^+$} ;
   \draw (2*2+4-8,4-2) node[scale=0.7] {$2\;6^+$} ;
  \foreach \y in {2,3} {
   \newcount\u ;
   \pgfmathsetcount{\u}{\x+\y} ;
   \draw[fl] (\y-8+2*\x,\y-2) -- (\y-7+2*\x,\y-1) ;
   \draw[fl] (\y-7+2*\x,\y-1) -- (\y-6+2*\x,\y-2) ;
  } ;
  \draw[fl] (-3,1) -- (-2,1) ;
  \draw (-2,1) node[scale=0.7] {$1\;5^-$} ;
  \draw[fl] (-2,1) -- (-1,1) ;
  \draw[fl] (-1,1) -- (-0.1,1) ;
  \draw (0,1) node[scale=0.7] {$2\;6^-$} ;
  \draw[fl] (0,1) -- (1,1) ;
} ;
\end{scope}
\begin{scope}[xshift=240]
\draw[thick, dashed, blue] (-5.2,-0.3) -- ++(3.8,0) -- ++(2.6,2.6) -- ++(-3.8,0) -- cycle ;
\end{scope}
\draw (-0.5,1) node[scale=1.5] {$\cdots$} ;
\draw (10,1) node[scale=1.5] {$\cdots$} ;
\end{tikzpicture}
\caption{The AR-quiver of $\cal D_{1,1}$}
\label{12}
\end{figure}
\begin{eqnarray*}
\X_{1}=\{[(0)] \}&\quad\quad& \X^\perp_{1}[-1]=\{\cal D_{1,1}\}\nonumber \\
\X_{2}=\{[(13)]\}&\quad\quad& \X^\perp_{2}[-1]=\{[(13)],[(15^{+})],[(15^{-})]\}\nonumber \\
\X_{3}=\{[(24)] \}&\quad\quad&\X^\perp_{3}[-1]=\{[(24)],[(26^{+})],[(26^{-})]\}\nonumber \\
\X_{4}=\{[(13)],[(15^{+})] \}&\quad\quad&\X^\perp_{4}[-1]=\{[(13)],[( 15^{-})]\}\nonumber \\
\X_{5}=\{[(24)],[(26^{+})] \}&\quad\quad& \X^\perp_{5}[-1]=\{[(24)],[( 26^{-})]\}\nonumber \\
\end{eqnarray*}
\end{example}
\vspace{-0.75cm}
Note that $\X_{4}=\{[(13)],[(15^{+})]\}$, $\X^\perp_{4}[-1]=\{[(13)],[( 15^{-})]\}$, $\X_{5}=\{[(24)],[(26^{+})] \}$, and $\X^\perp_{5}[-1]=\{[(24)],[( 26^{-})]\}$ are all the subcategories containing one diameter and some indecomposable rigid objects, its number is $4=2t_{1,1}$. All the torsion pairs $(\X_{i}, \X^\perp_{i})$ and $( ^\perp\X_{i}, \X_{i})$ are listed above for $i=1,2,3,4,5$.

\section{Hearts of torsion pairs}

In this section, we determine
the hearts of torsion pairs in finite $2$-CY triangulated categories with maximal rigid objects. Hearts of cotorsion pairs in any triangulated category were introduced by Nakaoka \cite{n1}, which unify
the construction of hearts of t-structures \cite{bbd} and construction of the abelian quotient categories
by cluster tilting subcategories \cite{bmrrt,kr,kz}. For two subcategories $\X,\Y$ in a triangulated category $\cal C$, the pair $(\X,\Y)$ is a torsion pair in $\cal C$ if and only if $(\X,\Y[-1])$ is a cotorsion pair in $\cal C$. The heart of torsion pair $(\X,\Y)$ is by definition the heart of cotorsion pair $(\X,\Y[-1])$. We will use the notation of cotorsion pairs in this section.

We recall the construction of hearts of cotorsion pairs from \cite{n1}: Let $\cal C$ be a triangulated category and $(\X,\Y)$ a cotorsion pair with core $\cal I$ in $\cal C$. Denote by $\cal H$ the subcategory
$(\cal X[-1]\ast \cal I)\cap (\cal I\ast\cal Y[1])$. The heart of the cotorsion pair $(\cal X,\cal Y)$ is
defined as the quotient category $\cal H/\cal I$, denoted by $\underline{\cal H}$.
\medskip

It was proved that $\underline{\cal H}$ is an abelian category \cite{n1}. There is a cohomological functor $H=h \pi$ from $\cal C$ to $\underline{\cal H}$, where $\pi$ is the quotient functor from $\cal C$ to $\underline{\cal C}=\cal C/\cal I$ and
$h$ is a functor from $\underline{\C}$ to $\underline{\H}$. See \cite{an,n1} for the details of the constructions. We give the  main result in this section.

\begin{proposition}
Let $\C$ be a finite $2$-Calabi-Yau triangulated category, and $(\X, \Y)$ be a cotorsion pair in $\C$ with core $\I=add I$, where $I$ is a rigid object, Then we have an equivalence of abelian categories
$$\underline{\cal H}\simeq mod \ End I$$
\label{g1}
\end{proposition}
\proof For any cotorsion pair $(\X,\Y)$  in $\C$ with core $\I$, $({}^\bot\I[1])/\I$ is also a finite $2$-CY triangulated category with shift functor $<1>$ \cite{iy}, and $(\X/\I,\Y/\I)$ is a t-structure by \cite{zz2}. It follows from Proposition \ref{p4} that ${}^\bot\I[1])/\I=\X/\I\bigoplus \Y/\I$.

On the other hand, $(\I,{}^\bot\I[1] )$ is a cotorsion pair with the same core $\I$, and the heart $\underline{H}_1$ of $(\I,{}^\bot\I[1] )$ is equivalent to the module category mod End $I$ by \cite{iy}, i.e, $\underline{H}_1\simeq \text{mod}\ EndI$ . By the same proof as Theorem 6.4 in \cite{zz2}, the heart $\underline{H}$ of $(\X,\Y)$ is equivalent to $\underline{H}_1$. Thus
$\underline{H}\simeq\mod\ EndI$. \qed

Now we have the following conclusions about the hearts of cotorsion pairs in finite $2$-CY triangulated categories $\cal C$:
1. If $\cal C$ contains cluster tilting objects, then hearts have been determined in \cite{zz2}. 2. If $\cal C$ has only zero maximal rigid objects, then any cotorsion pair $(\X,\Y)$ is a t-structure, and $\X[1]=\X,  \Y[1]=\Y.$ Then the heart $\H=\X[-1]\bigcap \Y[1]=0$. 3. If $\cal C$ has non-zero maximal rigid objects which are not cluster tilting, then the hearts are determined in the following result combining Proposition \ref{g1}.
\begin{corollary}
\begin{itemize}
  \item [1.]The heart of any cotorsion pair in $\A_{n,t}$ is module category over the algebras given by one of the following quivers with relations:
  \begin{itemize}
    \item [(1)]$\xymatrix@C=1.5em@R=1em{ 1 \ar@{->}[r] & 2 \ar@{->}[r] &  \ar@{..}[r]& k-1\ar@{->}[r]&k}$ with $1\leq k\leq n$.
    \item [(2)]$\xymatrix@C=1.5em@R=1em{ 1 \ar@{->}[r] & 2 \ar@{->}[r] &  \ar@{..}[r]& k-1 \ar@{->}[r] & k \ar@(ur,dr)^{\alpha} }$ with relation $\alpha^{2}$,  $1\leq k\leq n$.
    \item[(3)] Mutations of the quiver occurred in the above $(1)$ or $(2)$.

  \end{itemize}
  \item [2.] The heart of any cotorsion pair in $\D_{n,t}$ is module category over the algebras given by one of the following quivers with relations:
  \begin{itemize}
  \item [(1)]$\xymatrix@C=1.5em@R=1em{ 1 \ar@{->}[r] & 2 \ar@{->}[r] &  \ar@{..}[r]& k-1\ar@{->}[r]&k}$ with $1\leq k\leq n$.
    \item [(2)]$\xymatrix@C=1.5em@R=1em{ 1 \ar@{->}[r] & 2 \ar@{->}[r] &  \ar@{..}[r]& k-1 \ar@{->}[r] & k \ar@(ur,dr)^{\alpha} }$ with relation $\alpha^{2}$,  $1\leq k\leq n$.
    \item[(3)] Mutations of the quiver occurred in the above $(1)$ or $(2)$.
        \end{itemize}
   \item [3.] The heart of any cotorsion pair in $\D^{b}(\mathbb{K}E_{7})/\tau^{2}$ is module category over the algebras given by the following quiver with relation:
   $$\xymatrix@C=1.5em@R=1em{\cdot \ar@(ur,dr)^{\alpha} }$$ with relation $\alpha^{3}$.
   \item [4.] The heart of any cotorsion pair in $\D^{b}(\mathbb{K}E_{7})/\tau^{5}$ is module category over the algebras given by one of the following quivers with relations::
   \begin{itemize}
     \item [(1)] $\xymatrix@C=1.5em@R=1em{\cdot \ar@(ur,dr)^{\alpha} }$ with relation $\alpha^{2}$.
     \item [(2)]  $\xymatrix@C=1.5em@R=1em{\cdot \ar@(ul,dl)_{\alpha} \ar@{->}[r]^{\beta} & \cdot\ar@(ur,dr)^{\gamma} }$ with relations $\beta\alpha-\gamma\beta$, $\alpha^{2}$, $\gamma^{2}$.
   \end{itemize}
\end{itemize}
\end{corollary}
\proof By Proposition \ref{g1}, for any cotorsion pair $(\X,\Y)$ with core $\I$ in a finite $2$-CY triangulated category, the heart $\underline{H}\simeq mod \ End I$, $I$ is a rigid object in $\cal C$. So the heart is determined by the  endomorphism algebra of some rigid objects, a subalgebra of endomorphism algebra of a maximal rigid object. Since $[(1,3)]\bigoplus[(1,4)]\bigoplus\ldots\bigoplus[(1,n+2)]$ is a maximal rigid object in $\A_{n,t}$, and  its endomorphism algebra is given by $\xymatrix@C=1.5em@R=1em{ 1 \ar@{->}[r] & 2 \ar@{->}[r] & 3 \ar@{..}[r]& n-1 \ar@{->}[r] & n \ar@(ur,dr)^{\alpha} }$ with relation $\alpha^{2}$ by \cite{bpr}, we can get any endomorphism algebra of a maximal rigid object through some mutations, as a result, the endomorphism algebra of any rigid objects is obtained. We prove the assertion in $1.$

 Statement in $2.$ can be proved similarly. $\D^{b}(\mathbb{K}E_{7})/\tau^{2}$ has only two indecomposable rigid objects and each is maximal, the endomorphism algebra is given in Proposition $2.14$ in \cite{bpr}, so the result in $3.$ is clear. For $4.$, $\D^{b}(\mathbb{K}E_{7})/\tau^{5}$ has five indecomposable rigid objects, the endomorphism algebra of any maximal rigid object is given in Proposition $2.12$ in \cite{bpr}.\qed
 \bigskip

 $\bold {Acknowlegement}$: The first author expresses her great thankfulness to Panyue Zhou for his generous discussion on the topic of the paper. Both authors would like to thank Yu Zhou for his comments on the earlier version of the paper.

\end{document}